\documentclass[12pt]{article}
\usepackage{graphicx}
\usepackage{setspace}
\usepackage{bigints}
\usepackage{graphics} 
\usepackage{epstopdf}
\usepackage{epsfig} 
\usepackage{mathptmx} 
\usepackage{times} 
\usepackage{amsmath} 
\usepackage{amsthm}
\usepackage{amssymb}  
\usepackage{algorithmicx}
\usepackage{algorithm}
\usepackage{bm}
\usepackage[noend]{algpseudocode}
\usepackage{cite}
\usepackage{subfig}
\usepackage{gensymb}
\usepackage{booktabs}
\usepackage{url}

\usepackage{tikz,pgfplots}
\pgfplotsset{compat=newest} 
\pgfplotsset{plot coordinates/math parser=false} 

\DeclareMathOperator*{\minimize}{minimize}

\DeclareMathOperator*{\subjectto}{subject\:to}

\DeclareMathAlphabet{\pazocal}{OMS}{zplm}{m}{n}

\usepackage{amsfonts}
\newcommand{\overbar}[1]{\mkern 3.0mu\overline{\mkern-3.0mu#1\mkern-3.0mu}\mkern 3.0mu}

\newlength\figureheight 
\newlength\figurewidth  

\usepackage{algcompatible}
\usepackage{algpseudocode}

\makeatletter
\renewcommand*{\ALG@name}{Workflow}
\makeatother

\title{A predictive path-following controller for multi-steered articulated vehicles}

\author
{\small Oskar Ljungqvist and  Daniel Axehill\\
	\\
	\small{Department of Automatic Control, Link\"oping University, Link\"oping, Sweden}\\
	\small{E-mail: \texttt{\{oskar.ljungqvist, daniel.axehill\}@liu.se}}
}

\date{}

\begin{document}

	\baselineskip 16pt

	\maketitle 

\begin{abstract} 
Stabilizing multi-steered articulated vehicles in backward motion is a complex task for any human driver. Unless the vehicle is accurately steered, its structurally unstable joint-angle kinematics during reverse maneuvers can cause the vehicle segments to fold and enter a jack-knife state. In this work, a model predictive path-following controller is proposed enabling automatic low-speed steering control of multi-steered articulated vehicles, comprising a car-like tractor and an arbitrary number of trailers with passive or active steering. The proposed path-following controller is tailored to follow nominal paths that contains full state and control-input information, and is designed to satisfy various physical constraints on the vehicle states as well as saturations and rate limitations on the tractor's curvature and the trailer steering angles. The performance of the proposed model predictive path-following controller is evaluated in a set of simulations for a multi-steered 2-trailer with a car-like tractor where the last trailer has steerable wheels.
\end{abstract}

%
\section{Introduction}
The transportation sector faces growing demands from the society to increase efficiency and to reduce the environmental footprint related to freight and public transport. As a result, recent trends in modern transport include an increased interest in large capacity (multi-) articulated buses~\cite{michalek2019modular} and long tractor-trailer combinations~\cite{islam2015comparative}. In order to improve these long vehicle's maneuvering capability, some of the trailers are equipped with steerable wheels. Compared to single-steered N-trailer (SSNT) vehicles where all trailers are passive, multi-steered N-trailer (MSNT) vehicles are more agile, but also significantly more difficult to control for a human driver. This is partly because of the vehicle's additional degrees of freedom and due to specific kinematic and dynamics properties of MSNT vehicles~\cite{tilbury1995multisteering,islam2015comparative,orosco2002modeling,michalek2019modular}. To aid the driver, various control systems have been proposed to automatically control the steerable trailer wheels to either decrease the turning radius during reverse maneuvers or to diminish the so-called off-tracking effect during tight cornering~\cite{Beyersdorfer2013tractortrailer,odhams2011active,varga2018robust,van2015active,michalek2019modular}. 

Even though several feedback-control strategies have been proposed for various SSNT vehicles (see e.g.~\cite{michalek2014highly,LjungqvistJFR2019,hybridcontrol2001,rimmer2017implementation,altafini2003path}), only a limited amount of work has been devoted to the path-following or the trajectory-tracking control problem for special classes of MSNT vehicles (see e.g. \cite{odhams2011active,varga2018robust,van2015active,Yuan2015,sadeghi2019gain}). However, these approaches mainly use the additional trailer-steering capability to reduce the off-tracking effect while tracking a geometric reference path or trajectory. As a consequence, there is still a need to present a path-following controller for a generic MSNT vehicle for the case when the nominal path contains full state and control-input information, i.e., it is tailored to operate in series with a motion planner similar to~\cite{evestedtLjungqvist2016planning,li2019trajectory,Beyersdorfer2013tractortrailer,LjungqvistJFR2019}.   

The contribution of this work is a path-following controller for a generic MSNT with a car-like tractor targeting low-speed maneuvers, which is design to operate in series with a motion planner that computes feasible paths. It is done by first deriving a path-following error model describing the vehicle in terms of deviation from the nominal path. This error model together with physical constraints on states and control inputs are then used to design a path-following controller based on the framework of model predictive control (MPC)~\cite{mayne2000constrained,garcia1989model,faulwasser2015nonlinear,lima2017spatial}. To the best of the authors knowledge, this paper presents the first path-following controller for a generic MSNT with a car-like tractor admitting mixtures of off-axle/on-axle hitch connections and steerable/non-steerable trailers, and is designed to satisfy various constraints on states and control inputs.

The remainder of the paper is structured as follows. The kinematic vehicle model is presented in Section~\ref{c9:sec:model} and the path-following error model is derived in Section~\ref{c9:sec:error_model}. The proposed model predictive path-following controller is presented in Section~\ref{c9:sec:MPC}. In Section~\ref{c9:sec:results}, simulation results for a MS2T with a car-like tractor is presented and the paper is concluded in Section~\ref{c9:sec:conclusions} by summarizing the contributions and discussing directions for future work.
\section{Kinematic vehicle model}
\label{c9:sec:model}
\begin{figure}[t!]
	\centering
	\includegraphics[width=0.7\linewidth]{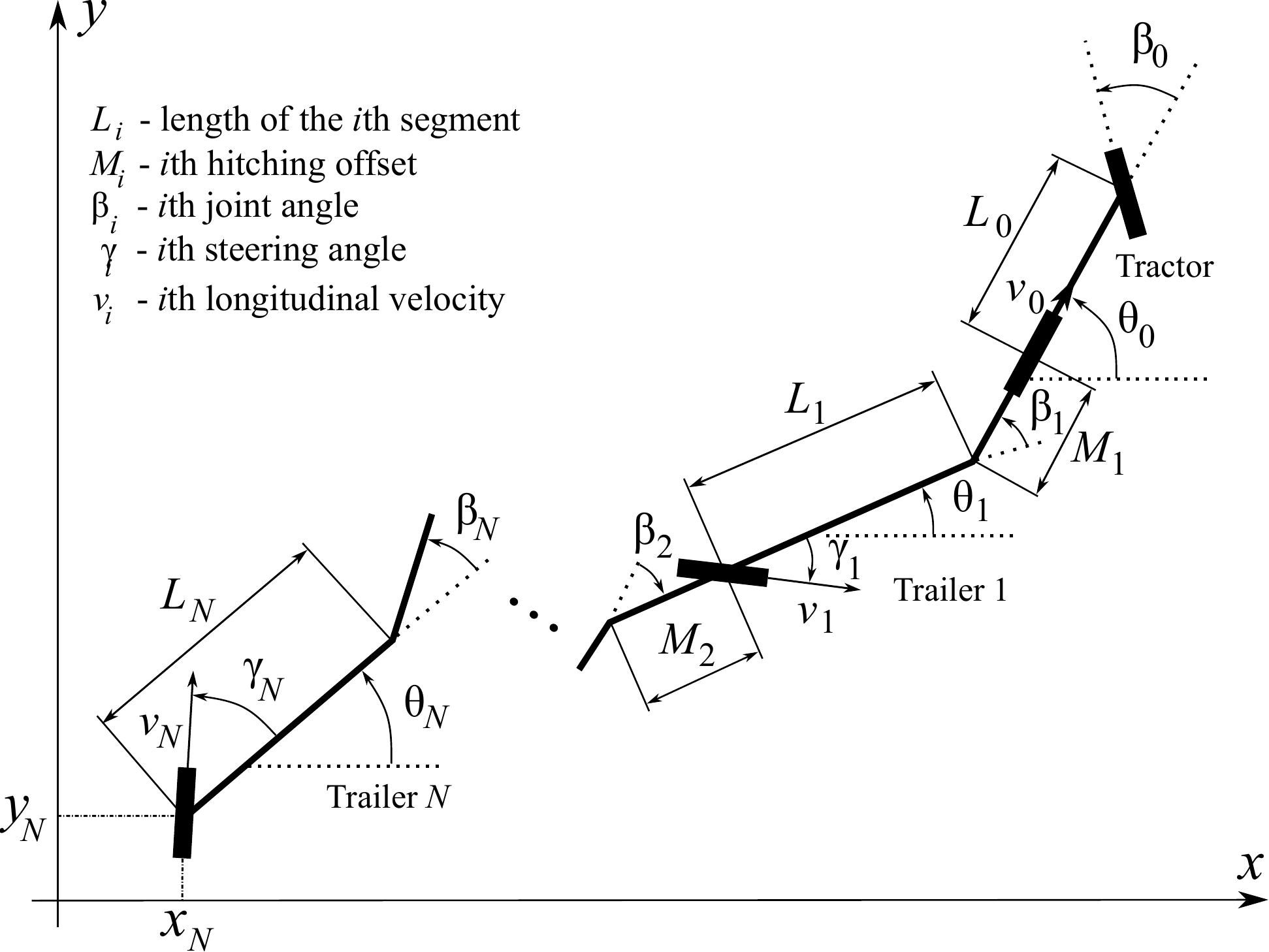}
	\caption{A schematic description of the MSNT with a car-like tractor in a global coordinate system.}
	\label{c9:fig:vehicle_model}%
\end{figure}

The MSNT with a car-like tractor considered in this work is composed of $N+1$ interconnected vehicle segments, including a leading car-like tractor and $N$ number of trailers that are either passively of actively steered. The car-like tractor has a steerable front wheel and its rear wheel is fixed. The MSNT vehicle is illustrated in Figure~\ref{c9:fig:vehicle_model}, where each vehicle segment is described by a segment length $L_i>0$ and a signed hitching offset $M_i$. Since low-speed maneuvers are considered, a kinematic model is used to describe the vehicle. The model is based on the work in~\cite{michalek2019modular} and is derived based on some assumptions including that wheels are rolling without slipping. By considering the steering angles as control inputs, the MSNT with a car-like tractor can be described with a state vector that consists of $n=3 + N$ variables:
\begin{itemize}
	\item[--] the global pose $(x_{N},y_{N},\theta_{N})$ of the $N$th trailer in a fixed coordinate frame
	\begin{align}
	q_{N} = \begin{bmatrix}
	  x_{N} & y_{N} & \theta_{N}
	\end{bmatrix}^T \in \mathbb R^2 \times \mathbb S,
	\label{c9:eq:poes_trailer}
	\end{align} 
	where $\mathbb S=(-\pi,\pi]$.
	\item[--] for $i=1,\hdots, N,$ a number of $N$ constrained joint angles
	\begin{align}
	\beta_i = \theta_{i-1} -\theta_{i} \in\pazocal B_i=[-\bar\beta_i,\bar\beta_i], \quad \bar \beta_i \in (0,\pi/2). 
	\label{c9:eq:joint_angles}
	\end{align}
	The state vector for the MSNT with a car-like tractor is defined as
	\begin{align}
	x = \begin{bmatrix}
	q_N^T & \beta_N & \beta_{N-1} & \hdots & \beta_{1}    
	\end{bmatrix}^T \in \pazocal X,
	\label{c9:eq:configuration_vector}
	\end{align}
	where 
	$\pazocal X =  \mathbb R^2 \times \mathbb S\times \pazocal B_N \times \pazocal B_{N-1} \times \hdots  \times\pazocal B_1.$
\end{itemize}
By treating the longitudinal velocity of the car-like tractor $v_0$ as an exogenous input, the control input consists of \mbox{$m=1+S$} variables:
\begin{itemize}
	\item[--] the curvature of the car-like tractor $\kappa_0 =\frac{\tan\beta_0}{L_0}$:
	\begin{align}
	\kappa_0 \in \pazocal Q_0=[-\bar\kappa_0,\bar\kappa_0], 
	\label{c9:eq:front_steering_angle}
	\end{align}
	where $\beta_0\in[-\bar\beta_0,\bar\beta_0]$, $\bar\beta_0\in(0,\pi/2)$ is the steering angle of the tractor's front wheels and \mbox{$\bar\kappa_0=\frac{\tan\bar\beta_0}{L_0}$} is the maximum curvature,
	\item[--] and $S\in\{1,\hdots,N\}$ number of steering angles associated with actively steered trailers 
	\begin{align}
	\gamma_a \in \pazocal Q_a=[-\bar\gamma_a,\bar\gamma_a], \quad \bar \gamma_a \in (0,\pi/2),
	\label{c9:eq:trailer_steering_angles}
	\end{align}
	where index $a\in\pazocal I_a \subseteq \{1,\hdots,N\}$ specifies which trailers that have steerable wheels. The control input for the MSNT with a car-like tractor is defined as
	\begin{align}
	u = \begin{bmatrix}
	\kappa_0 & \bm \gamma_a^T 
	\end{bmatrix}^T \in \pazocal U,
	\label{c9:eq:control_vector}
	\end{align}
	where $\kappa_0$ is the tractor's curvature and $\bm \gamma_a$ represents a vector of trailer steering angles and 
	$\pazocal U =\pazocal Q_0\times \underbrace{\pazocal Q_a\times \hdots\times \pazocal Q_a}_{S\text{-times}}$.
\end{itemize}
The leading car-like tractor is described by a kinematic single-track vehicle model and its orientation $\theta_0$ evolves as
\begin{align}
\dot \theta_0 = v_0\kappa_0.
\label{c9:eq:model_orientation_tractor}
\end{align}
Between any two neighboring vehicle segments, the transformation of the angular $\dot\theta_i$ and longitudinal $v_i$ velocities are given by (see, e.g.,~\cite{michalek2019modular}): 
\begin{align}
\begin{bmatrix}
\dot \theta_{i} \\  v_{i}
\end{bmatrix} &=
\underbrace{\begin{bmatrix}
	-\frac{M_{i}}{L_{i}}\frac{\cos{(\beta_{i}-\gamma_{i})}}{\cos\gamma_{i}} & \frac{\sin{(\beta_{i}-\gamma_{i}+\gamma_{i-1})}}{L_{i}\cos\gamma_{i}} \\[10pt]
	M_{i}\frac{\sin{\beta_{i}}}{\cos\gamma_{i}} & \frac{\cos{(\beta_{i}+\gamma_{i-1})}}{\cos\gamma_{i}}
	\end{bmatrix}}_{\triangleq J_i(\beta_{i},\gamma_i,\gamma_{i-1})}
\begin{bmatrix}
\dot \theta_{i-1} \\  v_{i-1}
\end{bmatrix},\quad i=1,\hdots,N,
\label{c9:eq:velocity_transformation}
\end{align}
where $\gamma_i$ denotes the steering angle of the $i$th trailer. Note that if the $j$th trailer is non-steerable, it suffices to take $\gamma_j=0$, and that $\gamma_0=0$ because the tractor's rear wheel is fixed. 

To satisfy actuator limitations, the rate of each trailer steering angle $\gamma_a$, $a\in\pazocal I_a$ and the tractor's curvature $\kappa_0$ are constrained as
\begin{align}
\begin{split}
\label{c9:eq:rate_constraints_on_steering}
|\dot\gamma_a| &\leq\dot{\bar\gamma}_a, \quad a\in\pazocal I_a, \\
|\dot\kappa_0| &\leq \dot{\bar\kappa}_{0},
\end{split}
\end{align}
which is compactly represented as $\dot { u}\in \Omega$.
Moreover, the position of the $N$th trailer evolves according to standard unicycle kinematics (see Figure~\ref{c9:fig:vehicle_model}) 
\begin{align}
\begin{split}
\label{c9:eq:position_Ntrailer}
\dot x_{N} &= v_{N}\cos(\theta_{N}+\gamma_N), \\
\dot y_{N} &= v_{N}\sin(\theta_{N}+\gamma_N). 
\end{split}
\end{align}
Using~\eqref{c9:eq:model_orientation_tractor} and~\eqref{c9:eq:velocity_transformation}, the angular rate $\dot \theta_{N}$ and longitudinal velocity $v_{N}$ of the $N$th trailer are given by
\begin{align}
\label{c9:eq:rate_speed_Ntrailer}
\begin{bmatrix}
\dot \theta_{N}\\
v_{N} 	
\end{bmatrix}=\prod_{i=0}^{N-1} \bm J_{N-i}(\beta_{N-i},\gamma_{N-i},\gamma_{N-i-1})\begin{bmatrix} v_0\kappa_0\\
v_0 
\end{bmatrix}.
\end{align}
Note that $v_0$ enters bilinearly in~\eqref{c9:eq:rate_speed_Ntrailer}. Therefore, using~\eqref{c9:eq:rate_speed_Ntrailer} and by introducing the vectors $ c^T=\begin{bmatrix}
1&0
\end{bmatrix}$ and $ d^T=\begin{bmatrix}
0&1
\end{bmatrix}$, the curvature of the $N$th trailer is defined as
\begin{align}
\label{c9:eq:curvature_Ntrailer}
\kappa_N(\beta_1,\hdots,\beta_N, u) \triangleq \frac{\dot\theta_N}{v_N} = 
\frac{ c^T\prod_{i=0}^{N-1} \bm J_{N-i}(\beta_{N-i},\gamma_{N-i},\gamma_{N-i-1})
	\begin{bmatrix} \kappa_0 \\ 1 \end{bmatrix}}
{ f_{v_N}(\beta_1,\hdots,\beta_N, u)},
\end{align}
where
\begin{align}
\label{c9:eq:velocity_0N_trans}
f_{v_N}(\beta_1,\hdots,\beta_N, u) =   d^T\prod_{i=0}^{N-1} \bm J_{N-i}(\beta_{N-i},\gamma_{N-i},\gamma_{N-i-1})
\begin{bmatrix} \kappa_0\\ 1\end{bmatrix},
\end{align}
which relates the longitudinal velocity transformation from the tractor $v_0$ to the $N$th trailer $v_N$ as $v_N=f_{v_N}(\beta_1,\hdots,\beta_N, u)v_0$. To guarantee that~\eqref{c9:eq:curvature_Ntrailer} is well defined, it is further assumed that the sets $\pazocal X$ and $\pazocal U$ are defined such that $f_{v_N}> 0$. Using~\eqref{c9:eq:position_Ntrailer}, \eqref{c9:eq:rate_speed_Ntrailer} and \eqref{c9:eq:velocity_0N_trans}, the model for the pose of the $N$th trailer can be represented as $\dot{ q}_N = v_Nf_{ q_N}( x, u)$.
Furthermore, using~\eqref{c9:eq:model_orientation_tractor} and \eqref{c9:eq:velocity_transformation}, the time derivative of~\eqref{c9:eq:joint_angles} yields the joint-angle kinematics
\begin{align}
\label{c9:eq:joint_angle_kinematics}
\dot\beta_i = &\dot\theta_{i-1} - \dot\theta_i =  c^T \prod_{j=N-i+1}^{N-1} \bm J_{N-j}(\beta_{N-j},\gamma_{N-j},\gamma_{N-j-1})\begin{bmatrix}
\kappa_0\\
1 
\end{bmatrix} \frac{v_N}{f_{v_N}(\beta_1,\hdots,\beta_N, u)}\nonumber \\
&-  c^T\prod_{j=N-i}^{N-1} \bm J_{N-j}(\beta_{N-j},\gamma_{N-j},\gamma_{N-j-1})\begin{bmatrix} \kappa_0\\
1 
\end{bmatrix}\frac{v_N}{f_{v_N}(\beta_1,\hdots,\beta_N, u)}, \quad i=1,\hdots,N.
\end{align}
Denote the joint-angle kinematics in~\eqref{c9:eq:joint_angle_kinematics} as $\dot\beta_i=v_Nf_{\beta_i}(\beta_1,\hdots,\beta_N, u)$, for $i=1,\hdots,N$. Now, the kinematic model of the MSNT vehicle with a car-like tractor is given in~\eqref{c9:eq:position_Ntrailer},~\eqref{c9:eq:rate_speed_Ntrailer} and \eqref{c9:eq:joint_angle_kinematics}, which can compactly be represented as
\begin{align}
\dot{ x}=v_N  f( x, u), \label{c9:vehicle_model}
\end{align}  
where $ f:\mathbb R^{n}\times \mathbb R^{m}\rightarrow \mathbb R^n$ is continuous and continuously differentiable with respect to $ x\in\pazocal X$ and $ u\in\pazocal U$. 

\section{Path-following error model}
\label{c9:sec:error_model}
\begin{figure}[t!]
	\centering
	\includegraphics[width=0.75\linewidth]{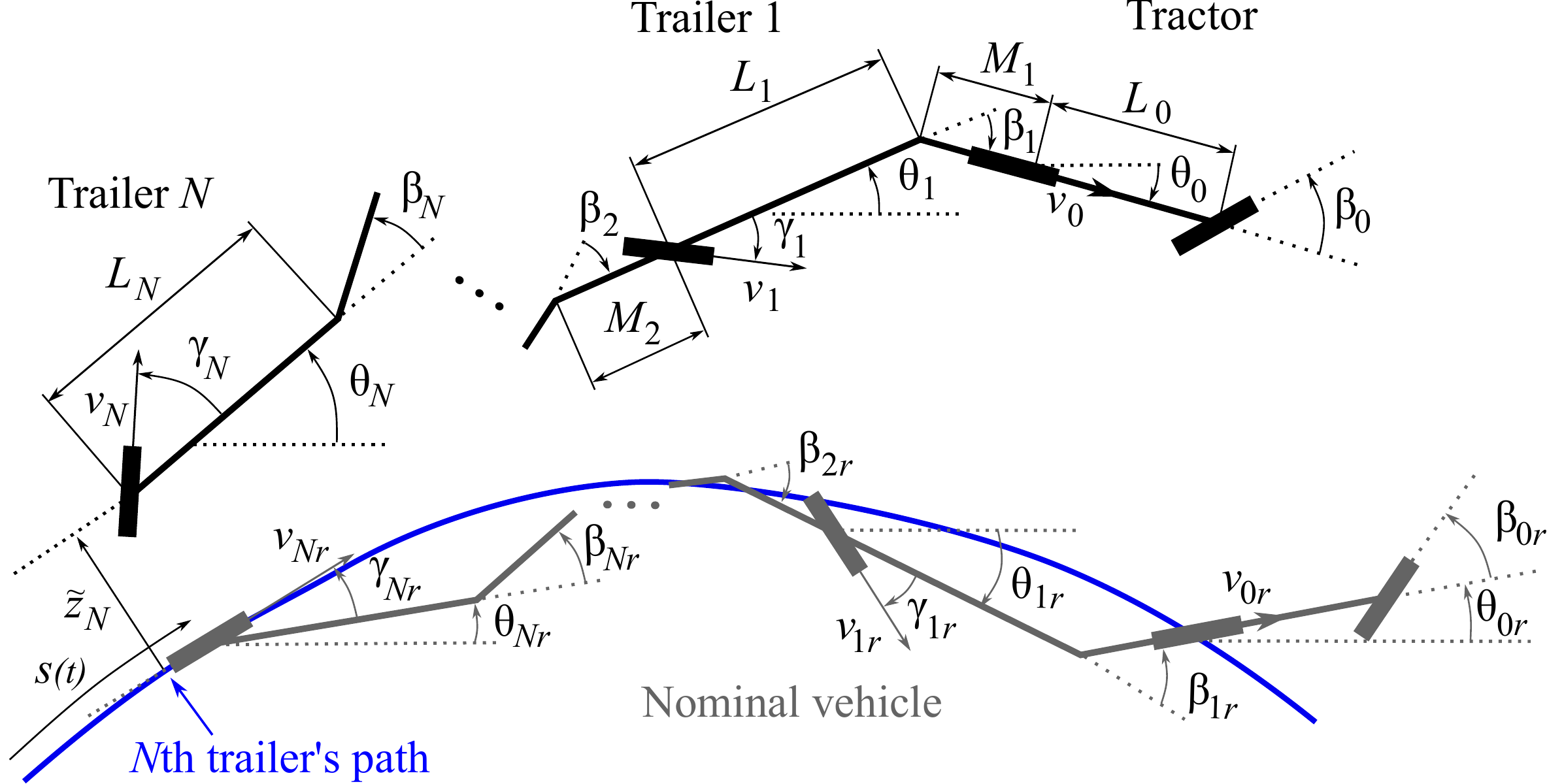}
	\caption{A schematic description of the MSNT with a car-like tractor in the Frenet frame coordinate system.}
	\label{c9:fig:vehicle_model_frenet_frame}
\end{figure}

It is assumed that a nominal trajectory $( x_r(\cdotp), u_r(\cdotp),v_{Nr}(\cdotp))$ for the MSNT vehicle~\eqref{c9:vehicle_model} is provided that satisfies the constraints on states ${ x}_r(\cdotp)\in {\pazocal X}$, and control inputs ${ u}_r(\cdotp)\in \pazocal U$ and \mbox{$\dot { u}_r(\cdotp)\in \Omega$}. Given the vehicle's current state $x(t)$, define $s(t)$ as the distance traveled by the position of the $N$th trailer onto its projection to its nominal path $(x_{Nr}(\cdotp),y_{Nr}(\cdotp))$ up to time $t$. By applying time-scaling~\cite{sampei1986time}, the nominal trajectory can instead be interpreted as a nominal path~\cite{LjungqvistJFR2019}:
\begin{align}\label{c9:nominal_path}
\frac{\text d  x_r}{\text d s} =  \bar v_{Nr}  f( x_r, u_r),
\end{align}
where $\bar v_{Nr}=\text{sign}{(v_{Nr})}\in\{-1,1\}$ specifies the nominal motion direction. Similar to~\cite{LjungqvistJFR2019}, the idea is now to model the MSNT vehicle in terms of deviation from this nominal path, as illustrated in Figure~\ref{c9:fig:vehicle_model_frenet_frame}. 
Denote $\tilde z_N(t)$ as the signed lateral error between the position of the $N$th trailer and its projection to its nominal path $(x_{Nr}(\cdotp), y_{Nr}(\cdotp))$. 
Denote the orientation error of the $N$th trailer as $\tilde{\theta}_N(t)=\theta_N(t)-\theta_{Nr}(s(t))$ and define the joint-angle errors $\tilde{\beta}_i(t)=\beta_i(t)-\beta_{ir}(s(t))$, $i=1,\hdots,N$. Finally, define the control-input deviation as \mbox{${\tilde{u}}(t) ={u}(t) - {u}_r(s(t))$} and let $\kappa_{Nr}=\kappa_{N}(\beta_{1r},\hdots,\beta_{Nr}, u_r)$ represent the curvature of the nominal path for the $N$th trailer.
Using the Frenet-frame transformation together with the chain rule, the MSNT vehicle~\eqref{c9:vehicle_model} can be described in terms of deviation from the nominal path~\eqref{c9:nominal_path} as 
\begin{subequations}%
	\label{c9:eq:model_frenet_frame}
	\begin{align}%
	\dot s = v_N& \frac{\bar v_{Nr}\cos (\tilde \theta_N+\tilde\gamma_N)}{1-\kappa_{Nr} \tilde z_N}, 
	\label{c9:eq:model_s1}
	\\ 
	\dot{\tilde z}_N = v_N& \sin(\tilde \theta_N+\tilde \gamma_N), 
	\label{c9:eq:model_s2}
	\\
	\dot{\tilde \theta}_N =  v_N&\Biggl(\kappa_N(\tilde\beta_1+\beta_{1r},\hdots,\tilde\beta_N+\beta_{Nr},{\tilde u}+ u_r) 
	- \frac{\kappa_{Nr}\cos (\tilde \theta_N+\tilde\gamma_N)}{1-\kappa_{Nr} \tilde z_N}\Biggr),
	\\
	\dot{\tilde \beta}_i = v_N&\Biggl(f_{\beta_i}(\tilde\beta_1+\beta_{1r},\hdots,\tilde\beta_N+\beta_{Nr},{\tilde u}+ u_r) \nonumber \\
	&-\frac{\cos (\tilde \theta_N+\tilde\gamma_N)}{1-\kappa_{Nr} \tilde z_N} f_{\beta_i}(\beta_{1r},\hdots,\beta_{Nr}, u_r)\Biggr), \quad i=N,N-1,\hdots,1.
	\label{c9:eq:model_s5}
	\end{align}
\end{subequations}
The transformation to the Frenet frame path-coordinate system is valid as long as $\tilde z_N$ and the sum $\tilde\theta_N+\tilde\gamma_N$ satisfy 
\begin{align}\label{c9:frenet_frame_transformation_constraints}
1-\kappa_{Nr}(s)\tilde z_N>0, \quad |\tilde\theta_N+\tilde\gamma_N|<\pi/2.
\end{align}
Essentially, this gives that $|\tilde z_N|<|\kappa_{Nr}^{-1}(s)|$ must hold when $\tilde z_N$ and $\kappa_{Nr}(s)$ have the same sign. Note that $\bar v_{Nr}$ is included in~\eqref{c9:eq:model_s1} to make $\dot s > 0$ as long as the constraints in~\eqref{c9:frenet_frame_transformation_constraints} are satisfied, and the $N$th trailer’s velocity $v_N$ and the nominal motion direction $\bar v_{Nr}$ have the same sign. Moreover, since it is assumed that $f_{v_N}> 0$  and the relationship $v_N = v_0f_{v_N}(\beta_1,\hdots,\beta_N, u)$ holds, an equivalent condition is that the velocity of the car-like tractor $v_0$ is selected such that $\text{sign}(v_0)=\bar v_{Nr}$.

Define the path-following error ${\tilde x} = [\tilde z_N \hspace{5pt} \tilde\theta_N \hspace{5pt} \tilde\beta_N \hspace{5pt}\hdots\hspace{5pt} \tilde\beta_1]^T$, where its model is given by~\eqref{c9:eq:model_s2}--\eqref{c9:eq:model_s5}. From the structure of~\eqref{c9:eq:model_s2}--\eqref{c9:eq:model_s5}, it is straightforward to verify that the origin $({\tilde x},{\tilde u})=(0,0)$ is an equilibrium point for all $t$. 
Moreover, since the velocity of the tractor $v_0$ is selected such that $\dot s(t)>0$, it is possible to perform time-scaling~\cite{sampei1986time} and eliminate the time-dependency presented in~\eqref{c9:eq:model_s2}--\eqref{c9:eq:model_s5}. Using the chain rule, it holds that $\frac{\text d\tilde x}{\text  ds}=\frac{\text d\tilde x}{\text  dt}\frac{1}{\dot s}$, and the spatial version of the path-following error model~\eqref{c9:eq:model_s2}--\eqref{c9:eq:model_s5} becomes 
\begin{subequations}%
	\label{c9:eq:model_frenet_frame_s}
	\begin{align}%
	\frac{\text d\tilde z_N}{\text d s} = \bar v_{Nr}& (1-\kappa_{Nr} \tilde z_N)\tan(\tilde \theta_N+\tilde \gamma_N), 
	\label{c9:eq:model_s5_spatia2}
	\\
	\frac{\text d \tilde \theta_N}{\text d s} =  \bar v_{Nr}& \biggl(\frac{1-\kappa_{Nr} \tilde z_N}{\cos (\tilde \theta_N+\tilde\gamma_N)}\kappa_N(\tilde\beta_1+\beta_{1r},\hdots,\tilde\beta_N+\beta_{Nr},{\tilde u}+ u_r) - \kappa_{Nr}\biggr),
	\\
	\frac{\text d\tilde \beta_i}{\text d s} = \bar v_{Nr}& \biggl(\frac{1-\kappa_{Nr} \tilde z_N}{\cos (\tilde \theta_N+\tilde\gamma_N)}f_{\beta_i}(\tilde\beta_1+\beta_{1r},\hdots,\tilde\beta_N+\beta_{Nr},{\tilde u}+ u_r) \nonumber \\ 
	&- f_{\beta_i}(\beta_{1r},\hdots,\beta_{Nr}, u_r)\biggr), \quad i=N,N-1,\hdots,1,
	\label{c9:eq:model_s5_spatial}
	\end{align}
\end{subequations}
which can be compactly represented as
\begin{align}
\frac{\text d {\tilde x} }{\text d s}=\bar v_{Nr}\tilde f(s,{\tilde x}, {\tilde u}),
\label{c9:eq:path_following_error_model}
\end{align}
where $\tilde f(s,0,0)=0$ for all $s$. In the next section, a model predictive path-following controller is proposed to stabilize the path-following error model~\eqref{c9:eq:path_following_error_model} around the origin, i.e., around the nominal path~\eqref{c9:nominal_path}. 


\section{Model predictive path-following controller}
\label{c9:sec:MPC}
The objective of the model predictive path-following controller is to control the tractor's curvature $\kappa_{0}$ and the trailer steering angles $\bm\gamma_a$ such that the path-following error is minimized,
while the constraints on states $x\in\pazocal X$, and control inputs $u\in\pazocal U$ and $\dot u\in\Omega$ are satisfied for all time instances. To obtain an MPC problem that can be solved online at a high sampling rate, the
goal is to derive an MPC formulation that can be
converted into the form of a quadratic programming (QP)
problem. First, the nonlinear path-following error model~\eqref{c9:eq:path_following_error_model}
is linearized around the origin $({\tilde x},{\tilde u}) = (0, 0)$, i.e., around the nominal path:
\begin{align}
\label{c9:linear_cont_model}
\frac{\text d {\tilde x}}{\text ds} = \bar v_{Nr} A(s){\tilde x} + \bar v_{Nr} B(s){\tilde u},
\end{align} 
where ${\tilde x}$ is the path-following error and $\tilde u$ is the control-input deviation. Using Euler-forward discretization with sampling distance $\Delta_s$, the discrete-time approximation of~\eqref{c9:linear_cont_model} becomes
\begin{align}
\label{c9:d_lin_sys}
{\tilde x}_{k+1} =  F_k{\tilde x}_k +  G_k{\tilde u}_k,
\end{align}  
where
\begin{align}
\label{c9:d_system_matrix}
 F_k =  I +\Delta_s \bar v_{Nr} A_k, \quad
 G_k = \Delta_s\bar v_{Nr} B_k.
\end{align} 
Since the tractor's curvature and the trailer-steering angles are $ u_k={\tilde u}_k+{ u}_{r,k}\in\pazocal U$, they are bounded as
\begin{align}
\label{c9:curvature_constaint_MPC}
-{\bar u}_k &\leq {\tilde u}_{k} + {u}_{r,k}  \leq {\bar u}_k, 
\end{align}
where ${\bar u}_k = [\bar \kappa_{0,k}\hspace{5pt} \hspace{5pt} \bm{\bar \gamma}^T_{a,k}]^T$ and $\bm{\bar \gamma}_{a,k}$ represents a vector of maximum trailer-steering angles.  
Furthermore, since $\dot s>0$ and $v_N = v_0 f_{v_N}(\beta_1,\hdots,\beta_N, u)$, the rate limits on the control input $\dot{ u}\in\Omega$ can be described in $s$ using the chain rule as 
\begin{align}
\begin{split}
\label{c9:rate-limit-constraints-ds}
\left|\frac{\text d\kappa_0}{\text ds}\right|&\leq\frac{\dot{\bar\kappa}_{0}}{\dot s} = \frac{1-\kappa_{Nr} \tilde z_N}{|v_0| f_{v_N}(\beta_1,\hdots,\beta_N, u)\cos(\tilde\theta_N+\tilde\gamma_N)}\dot {\bar\kappa}_{0}, \\
\left|\frac{\text d\gamma_a}{\text ds}\right|&\leq\frac{\dot{\bar\gamma}_{a}}{\dot s} = \frac{1-\kappa_{Nr} \tilde z_N}{|v_0| f_{v_N}(\beta_1,\hdots,\beta_N, u)\cos(\tilde\theta_N+\tilde\gamma_N)}\dot {\bar\gamma}_{a},
\end{split}
\end{align}
for $a\in\pazocal I_a$. Locally around the origin $({\tilde x},{\tilde u})=(0,0)$, it holds that $\cos(\tilde\theta_N+\tilde\gamma_N)\approx 1$ and $\kappa_{Nr} \tilde z_N\approx 0$. Therefore, to avoid coupling between $\tilde u$ and $\tilde x$, the constraints in~\eqref{c9:rate-limit-constraints-ds} are approximated as
\begin{align}
\begin{split}
\label{c9:rate-limit-constraints-ds2}
\left|\frac{\text d\kappa_0}{\text ds}\right| &\leq \frac{\dot {\bar\kappa}_{0}}{|v_0|f_{v_N}(\beta_{1r},\hdots,\beta_{Nr}, u_r)} \triangleq \bar c_{0}(s),\\
\left|\frac{\text d\gamma_a}{\text ds}\right| &\leq \frac{\dot {\bar\gamma}_{a}}{|v_0|f_{v_N}(\beta_{1r},\hdots,\beta_{Nr}, u_r)} \triangleq \bar c_{a}(s),\quad a\in\pazocal I_a.
\end{split}
\end{align} 
By discretizing~\eqref{c9:rate-limit-constraints-ds2} using Euler forward with sampling distance $\Delta_s$, the rate limits on the control input can be described by the following slew-rate constraint
\begin{align}
\label{c9:curvature_rate_constaint_MPC}
-\mathbf{\bar c}_{k}\Delta_s \leq {\tilde u}_{k} - {\tilde u}_{k-1} - {\bar u}_{r,k} \leq  \mathbf{\bar c}_k\Delta_s,   
\end{align}
where $\bar u_{r,k} = u_{r,k} - u_{r,k-1}$ and $\mathbf{\bar c}_{k} = [\bar c_{0,k}\hspace{5pt} \hspace{5pt} \mathbf{\bar c}^T_{a,k}]^T$, where $\mathbf{\bar c}_{a,k}$ represents a vector of rate limits for the trailer-steering angles. Denote the linear inequality constraints in~\eqref{c9:curvature_constaint_MPC} and~\eqref{c9:curvature_rate_constaint_MPC} as \mbox{${\tilde u}_k\in\tilde {\pazocal U}_k$}. 
Finally, since \mbox{$\beta_{i,k}=\beta_{ir,k}+\tilde\beta_{i,k}$}, $i=1,\hdots,N$, the constraints on the joint angles can be written as
\begin{align}
\label{c9:joint_angle_constaints_MPC}
-\bar\beta_i \leq \beta_{ir,k}+\tilde\beta_{i,k} \leq  \bar\beta_i, \quad i=1,\hdots,N,
\end{align}
which is compactly denoted as ${\tilde x}_k\in\tilde{\pazocal X}_k$. Note that $\tilde{\pazocal X}_k$ can be designed to also include constraints on other path-following error states. Now, given the path-following error $\tilde x(s(t))$ at time $t$, the MPC problem with prediction horizon $\overbar N$ is defined as follows
\begin{align}
\label{c9:MPC_problem}
\begin{split}
\minimize_{\mathbf{\tilde x},\hspace{0.1ex} \mathbf{\tilde u}} \hspace{3ex} &V_{\overbar N}(\mathbf{\tilde x},\mathbf{\tilde u}) = V_f(\tilde x_{\overbar N}) + \sum_{k=0}^{{\overbar N}-1} l(\tilde x_k,\tilde u_k) \\
\subjectto \hspace{2ex} 
&\tilde x_{k+1} =  F_k\tilde x_k +  G_k\tilde u_k, \hspace{2.5pt} k = 0,\hdots,{\overbar N}-1,\\
&{\tilde x}_k\in\tilde{\pazocal X}_k,  
\quad \tilde u_k\in\tilde{\pazocal U}_k,\hspace{10pt} k = 0,\hdots,{\overbar N}-1, \\
&\tilde x_0  = \tilde x(s(t))  \text{ given,}
\end{split}
\end{align}
where $\mathbf{\tilde x}^T = [\tilde x_0^T \hspace{5pt} \tilde x_1^T \hspace{5pt} \hspace{5pt}\hdots \hspace{5pt} \tilde x_{\overbar N}^T]$ is the path-following error sequence and $\mathbf{\tilde u}^T = [\tilde u_0^T \hspace{5pt} \tilde u_1^T \hspace{5pt} \hdots \hspace{5pt} \tilde u_{{\overbar N}-1}^T]$ is the input deviation sequence. The stage-cost is chosen to be quadratic $l(\tilde x_k,\tilde u_k) = ||\tilde x_k||^2_{Q} + ||\tilde u_k||_{R}^2$ as well as and the terminal cost $V_f(\tilde x_{\overbar N})=\tilde x_{\overbar N}^T P_{\overbar N}\tilde x_{\overbar N}$, where the matrices $ Q\succeq 0$, $ R\succ 0$ and $P_{\overbar N}\succ 0$ are design choices. Since the cost function $V_{\overbar N}$ is quadratic and there are only linear equality and inequality constraints, the optimization problem in~\eqref{c9:MPC_problem} can be written as a standard QP problem. Thus, at each sampling instance, the QP problem in~\eqref{c9:MPC_problem} is solved to obtain the optimal open-loop control-input deviation sequence $\mathbf{\tilde u}^*$. Only the first control-input deviation $\tilde u_0^*$ is deployed to the vehicle
\begin{align}
u(t) = u_r(s(t)) + \tilde u_0^*, 
\end{align}
and the QP problem~\eqref{c9:MPC_problem} is repeatedly solved at a fixed controller frequency $f_s$ using the current state estimate. Note that the MPC controller only computes the feedback part of the control input $\tilde u_0^*$, as the optimal feedforward $u_r(s(t))$ already is provided by the motion planner.
\subsection{Controller design}
We now turn to the problem of designing the cost function $V_{\overbar N}$ for the MPC controller~\eqref{c9:MPC_problem}. Since the nominal path contains full state and control-input information, it is possible to compute the nominal path as well as the nominal orientation of each vehicle segment using holonomic relationships~\cite{altafini1998general}. In order to minimize the risk of colliding with any obstacle, it is preferred that the MPC controller is tuned such that all path-following errors are penalized. Denote $\tilde z_i$, $i=0,\hdots,N-1$ as the lateral error of the $i$th vehicle segment with respect to its nominal path, and denote $\tilde \theta_i=\theta_i-\theta_{ir}$, $i=0,\hdots,N-1$ as their corresponding heading errors. As explained in~\cite{altafini2003path}, it is for general paths not possible to derive closed-form expressions to relate these auxiliary path-following errors as a function of the modeled ones $\tilde x$. However, around a straight nominal paths, closed-form expressions exist and the lateral and heading errors can be described as a function of $\tilde x$ using the following recursion 
\begin{align}
\begin{split}
\tilde z_i &= \tilde z_{i+1} + L_{i+1}\sin\tilde\theta_{i+1} + M_{i+1}\sin(\tilde\theta_{i+1}+\tilde\beta_{i+1}), \\
\tilde \theta_i &= \tilde\theta_{i+1} + \tilde\beta_{i+1}, \quad i = N-1,\hdots,0.
\end{split}
\end{align}   
Using these approximate relationships also for curved nominal paths, define
\begin{align}
\mathbf z = [\tilde x^T \hspace{5pt} \tilde z_{N-1} \hspace{5pt} \tilde \theta_{N-1} \hspace{5pt}\hdots \hspace{5pt} \tilde z_{0} \hspace{5pt} \tilde \theta_{0}]^T \triangleq  h_z(\tilde x),
\end{align}  
where $h_z(0)=0$, which defines the control-objective vector intended to be penalized. Since $ h_z(\tilde x)$ in nonlinear, it is linearized around the origin which yields \mbox{$\mathbf z = \frac{\partial  h_z(0)}{\partial \tilde x}\tilde x \triangleq  M\tilde x$}. The matrix $ M$ is then used to select the weight matrix for the quadratic stage-cost on $\tilde x$ as $ Q= M^T\bar Q  M$, where $\bar Q\succeq 0$ is a diagonal design matrix. Now, each diagonal element in $\bar Q$ penalizes a specific control objective in $ z$. The matrix $ M$ is then used to transform the specified design choice to $Q$, which typically obtains nonzero off-diagonal elements. 

When the matrices $ Q$ and $ R$ are selected, the weight matrix for the terminal cost $ P_{\overbar N}\succ 0$ is computed by solving the discrete-time algebraic Riccati equation (DARE):
\begin{align}\label{c9:DARE}
 P_{\overbar N} =  F^T P_{\overbar N}  F + Q -  F^T P_{\overbar N}  G K,
\end{align}
where $K =(R+ G^T P_{\overbar N} G)^{-1} G^T P_{\overbar N} F$ is the linear quadratic (LQ) feedback gain, and $F$ and $G$ are the discrete system matrices~\eqref{c9:d_system_matrix} for the linearized path-following error model~\eqref{c9:linear_cont_model} around a straight nominal path. Note that since the nominal motion direction \mbox{$\bar v_{Nr}\in\{-1,1\}$} enters bilinearly in~\eqref{c9:linear_cont_model}, the system's stability properties depend on the nominal motion direction. As a consequence, different terminal costs are used during backward and forward motion tasks~\cite{LjungqvistJFR2019,hybridcontrol2001}. Moreover, since the prediction model used in the MPC controller is an approximation, the originally hard joint-angle constraints are replaced with soft constraints using standard techniques~\cite{mayne2000constrained}.

Even though there exists a well-established theory for guaranteeing closed-loop stability for MPC (see e.g.~\cite{mayne2000constrained,garcia1989model}), a formal stability analysis is out of the scope in this work. Instead, extensive simulation trails are included to indicate the performance and stability properties of the MSNT vehicle using the proposed MPC controller~\eqref{c9:MPC_problem}.

\section{Simulation results}
\label{c9:sec:results}
In this section, the proposed model predictive path-following controller is evaluated on a MS2T with a car-like tractor, where trailer $N=2$ is steerable, i.e., $\pazocal I_a = \{2\}$, and a mixture of off-axle ($M_1\neq 0$) and on-axle ($M_2=0$) hitch connections. The vehicle parameters are presented in Table~\ref{c9:tab:vehicle_parameters}. Except for that trailer 2 is steerable, the parameters coincide with the full-scale test platform presented in~\cite{LjungqvistJFR2019}. Using the formulas presented in Section~\ref{c9:sec:model}, it is now straightforward to derive the kinematic vehicle model~\eqref{c9:vehicle_model} with state vector \mbox{$x=[x_2\hspace{5pt}y_2\hspace{5pt}\theta_2\hspace{5pt}\beta_2\hspace{5pt}\beta_1]^T$} and control input $u = [\kappa_0 \hspace{5pt} \gamma_2]^T$. Moreover, the path-following error is $\tilde x=[\tilde z_2 \hspace{5pt}\tilde \theta_2 \hspace{5pt} \tilde\beta_2\hspace{5pt} \tilde \beta_1]^T$ and the control-input deviation is $\tilde u=[\tilde \kappa_0\hspace{5pt}\tilde \gamma_2]^T$. Using the recursive formulas presented in Section~\ref{c9:sec:model} and Section~\ref{c9:sec:error_model} for this specific MS2T vehicle, the spatial path-following error model~\eqref{c9:eq:path_following_error_model} becomes
\begin{align}%
	\label{c9:eq:model_frenet_frame_2Ts}
	\begin{split}%
	\frac{\text d\tilde z_2}{\text d s} = \bar v_{2r} &(1-\kappa_{2r} \tilde z_2)\tan(\tilde \theta_2+\tilde \gamma_2), 
	\\
	\frac{\text d \tilde \theta_2}{\text d s} =  \bar v_{2r}&\left(\frac{1-\kappa_{2r} \tilde z_2}{\cos (\tilde \theta_2+\tilde\gamma_2)}\kappa_2(\beta_{2r}+\tilde\beta_{2},\gamma_{2r}+\tilde\gamma_2) - \kappa_{2r}\right),
	\\
	\frac{\text d\tilde \beta_2}{\text d s} = \bar v_{2r} &\Biggl(\frac{1-\kappa_{2r} \tilde z_2}{\cos (\tilde \theta_2+\tilde\gamma_2)}f_{\beta_2}(\beta_{1r}+\tilde\beta_1,\beta_{2r}+\tilde\beta_2,u_{r}+\tilde u)  
	-  f_{\beta_2}(\beta_{1r},\beta_{2r},u_r)\Biggr),
	\\
	\frac{\text d\tilde \beta_1}{\text d s} = \bar v_{2r} &\Biggl(\frac{1-\kappa_{2r} \tilde z_2}{\cos (\tilde \theta_2+\tilde\gamma_2)}f_{\beta_1}(\beta_{1r}+\tilde\beta_1,\beta_{2r}+\tilde\beta_2,u_{r}+\tilde u) 
	-  f_{\beta_1}(\beta_{1r},\beta_{2r},u_r)\Biggr),
	\end{split}
\end{align}
where $\kappa_{2r}=\kappa_2(\beta_{2r},\gamma_{2r})$ is the nominal curvature of trailer 2, and the functions $f_{\beta_1}$, $f_{\beta_2}$ and $\kappa_2$ are provided in Appendix A. The model in~\eqref{c9:eq:model_frenet_frame_2Ts} can compactly be written as $\text d {\tilde x}/\text d s = \bar v_{2r}\tilde f(s,\tilde x,\tilde u)$, and its linearization around the origin $(\tilde x,\tilde u)=(0,0)$ can be written as in~\eqref{c9:linear_cont_model}, where the matrices $A(s)$ and $B(s)$ are provided in Appendix A. The linearized system~\eqref{c9:linear_cont_model} is then discretized using a sampling distance $\Delta_s = 0.2$ m to obtain a discrete-time representation~\eqref{c9:d_lin_sys}. 

\begin{table}[t!]
	\caption{Vehicle parameters for the MS2T vehicle.}
	\centering
	\begin{tabular}{l l}
		\hline \noalign{\smallskip} Vehicle parameter  & Value   \\  \hline \noalign{\smallskip}	
		Tractor's wheelbase $L_0$        &   4.62 m  \\ 
		Length of off-hitch $M_1$        &   1.66 m  \\
		Length of trailer 1 $L_1$            &   3.87 m  \\ 
		Length of trailer 2 $L_2$            &   8.0 m  \\
		Maximum joint angles $\bar\beta_i$, $i=1,2$ & 0.8 rad \\
		Maximum curvature of tractor $\bar\kappa_0$ & 0.18 m$^{-1}$ \\   
		Maximum curvature rate of tractor $\dot{\bar\kappa}_0$ & $0.13$ $\text{m}^{-1}\text{s}^{-1}$\\   
		Maximum steering angle trailer 2  $\bar\gamma_3$ & 0.35 rad \\   
		Maximum steering-angle rate trailer 2 $\dot{\bar\gamma}_2$ & $0.8$ rad/s\\   
		\hline \noalign{\smallskip}
	\end{tabular}
	\label{c9:tab:vehicle_parameters}
\end{table}

\begin{table}[t!]
	\caption{Design parameters for the MPC controller.}
	\centering
	\begin{tabular}{l l}
		\hline \noalign{\smallskip} Vehicle parameter  & Value   \\  \hline \noalign{\smallskip}	
		Prediction horizon ${\overbar N}$       & 40  \\ 
		Weight matrix $\bar Q$  & $1/35\times \text{diag}([ 0.5\hspace{5pt} 1\hspace{5pt} 4 
		\hspace{5pt} 4 \hspace{5pt} 0.5 \hspace{5pt}1 \hspace{5pt}  0.5\hspace{5pt} 1])$  \\
		Weight matrix $R$       & $\text{diag}([4\hspace{5pt} 3])$  \\
		Sampling distance $\Delta_s$ & 0.2 m  \\ 
		Controller frequency $f_s$   & 10 Hz  \\
		\hline \noalign{\smallskip}
	\end{tabular}
	\label{c9:tab:MPC_parameters}
\end{table}

The proposed MPC controller is designed following the approach presented in Section~\ref{c9:sec:MPC}, where design parameters are in Table~\ref{c9:tab:MPC_parameters} and the control objective is $\mathbf z=[\tilde x^T\hspace{5pt}\tilde z_1\hspace{5pt}\tilde\theta_1\hspace{5pt}\tilde z_0\hspace{5pt}\tilde\theta_0]^T$. The terminal costs $P_{\overbar N}$ (one for forward and one for backward motion tasks) are computed by solving the DARE in~\eqref{c9:DARE} using the discrete system matrices $F= I + \Delta_s\bar v_{2r} A$ and $ G=\Delta_s\bar v_{2r} B$, obtained around a straight nominal path in forward ($\bar v_{2r}=1$) and backward ($\bar v_{2r}=-1$) motion. The matrices $A$ and $B$ for this special case are provided in Appendix A.  

\begin{figure}[b!]
	\centering
	\includegraphics[width=0.95\linewidth]{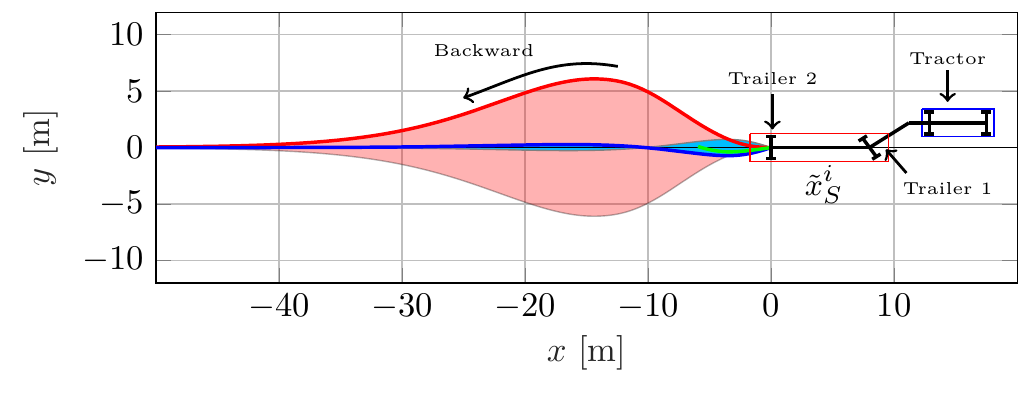}
	\vspace{-15pt}
	\caption{Path following of a straight nominal path ($y_{2r}=0$, black line) in backward motion from perturbed initial joint-angle errors $\tilde\beta^i_1,\tilde\beta^i_2\in[-0.6,\hspace{1pt}0.6]$ rad. The blue and red sets represent the convex envelope of the trajectories for the position of \mbox{trailer 2} using MS2T-MPC and SS2T-MPC, respectively. For the high-lighted initial state $\tilde x^i_{S}$, the paths taken by the position of trailer 2 are plotted for MS2T-LQ (green line), SS2T-MPC (red line) and  MS2T-MPC (blue line). As can be seen, MS2T-LQ leads to jack-knife (see also Figure~\ref{c9:fig:b23_traj_straight}).}
	\label{c9:fig:xy_path_straight}
\end{figure}

The MPC controller is implemented in Matlab using YALMIP where Gurobi 8.1.1 is used as QP solver~\cite{gurobi} to solve~\eqref{c9:MPC_problem} at each sampling instance. The performance of the proposed MPC controller is evaluated in a simulation study containing a straight and a figure-eight nominal path. The simulations are performed on a standard laptop computer with an Intel Core i7-4600U@2.1GHz CPU. The proposed MPC controller (MS2T-MPC) is benchmarked with an LQ controller (MS2T-LQ), as proposed in~\cite{LjungqvistJFR2019}. The LQ controller is given by $\tilde u=K\tilde x$, where the feedback gain $K$ is computed by solving the DARE in~\eqref{c9:DARE} using the weight matrices $Q$ and $ R$ that are also used by the MPC controller. Additionally, to analyze if the MPC controller is able to exploit the additional steering capability, it is also compared with an MPC controller for an SS2T vehicle, i.e., $\bar\gamma_2=0$, with the same vehicle parameters. This MPC controller (SS2T-MPC) uses the same design parameters except that the weight matrix on the control-input deviation is selected as $R=4$, because the control input $\tilde u=\tilde \kappa_0$ is a scalar for SS2T. Moreover, to make fair comparisons the nominal paths are designed to be feasible for the SS2T vehicle, i.e., $\gamma_{2r}=0$. In the simulations, the initial path-following error $\tilde x^i=\tilde x(0)$ is perturbed to compare how the different controllers handle disturbance rejection while satisfying the constraints on the joint angles.

\begin{figure}[t!]
	\centering
	\includegraphics[width=0.55\linewidth]{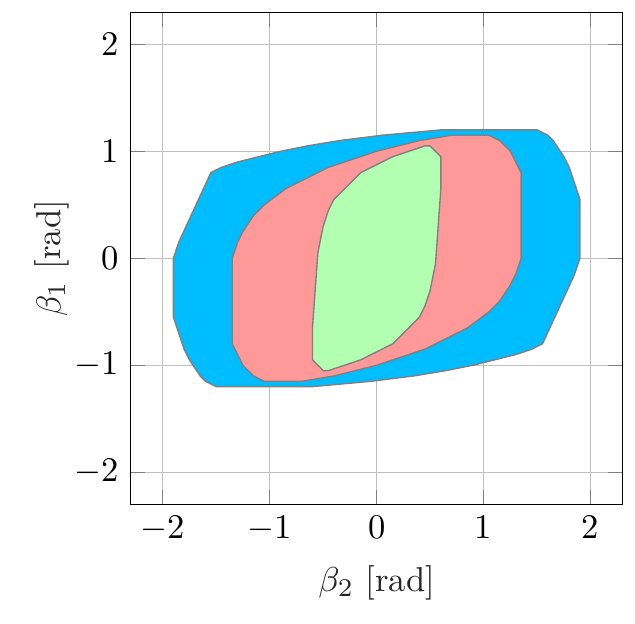}
	\vspace{-10pt}
	\caption{The simulated region of attraction for MS2T-LQ (green set), SS2T-MPC (red set) and MS2T-MPC (blue set) while following a straight nominal path ($y_{2r}=0$) in backward motion from perturbed initial joint-angles.}
	\label{c9:fig:b23_stability_region}
	\vspace{-0pt}
\end{figure}

The first set of simulations involves backward tracking of a straight nominal path aligned with the $x$-axis, where the longitudinal velocity of the tractor is selected as $v_0 = -1$ m/s. In this scenario, the initial state is \mbox{$\tilde x^i=[0\hspace{5pt}0\hspace{5pt}\beta_{2}^i\hspace{5pt}\beta^i_{1}]^T$}, where the initial joint-angle errors $\beta_{2}^i$ and $\beta_{1}^i$ are perturbed to various degrees. The setup is illustrated in Figure~\ref{c9:fig:xy_path_straight}. 
First, to analyze the stability region of the closed-loop systems, we numerically compute the region of attraction for the systems by performing simulations from a large set of initial joint-angle errors. In these simulations, the joint-angle constraints are temporarily removed from the MPC controllers and the closed-loop systems are checked for convergence to the straight nominal path. The resulting regions obtained from simulations are illustrated in Figure~\ref{c9:fig:b23_stability_region}. As expected, MS2T-MPC (blue set) has the largest region. Even though SS2T-MPC has non-steerable trailers, its region of attraction (red set) is larger than for MS2T-LQ (green set). This result is obtained because the LQ controller is not aware of the control-input constraints, as opposed to the two MPC controllers.

\begin{figure}[t!]
	\centering
	\captionsetup[subfloat]{captionskip=-3pt} 
	\setlength\figureheight{0.2\columnwidth}
	\setlength\figurewidth{0.33\columnwidth}
	\subfloat[][The joint-angle trajectories from initial state $\tilde x^i_{S}$ in Figure~\ref{c9:fig:xy_path_straight}. Initial (desired) state denoted by a black (blue) star.]{
		\begin{tikzpicture}
		\node[anchor=south west] (myplot) at (0,0) {
			\input{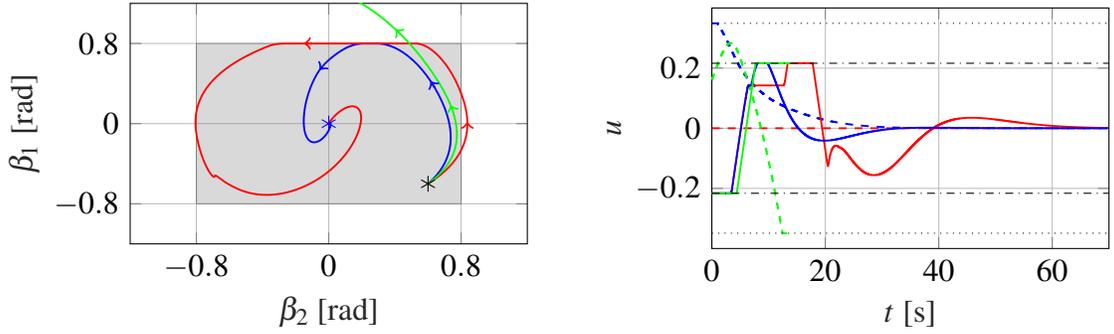}
		};
		\end{tikzpicture}
		\label{c9:fig:b23_traj_straight}
	}
	~
	\setlength\figureheight{0.2\textwidth}
	\setlength\figurewidth{0.33\textwidth}
	\subfloat[][Control inputs $\kappa_0$ (solid) and $\gamma_2$ (dashed) from ${\tilde x}^i_{S}$, and their limits $\bar\kappa_0$ (dashed-dotted black) and $\bar\gamma_2$ (dotted black) .]{
		\begin{tikzpicture}
		\node[anchor=south west] (myplot) at (0,0) {
			\input{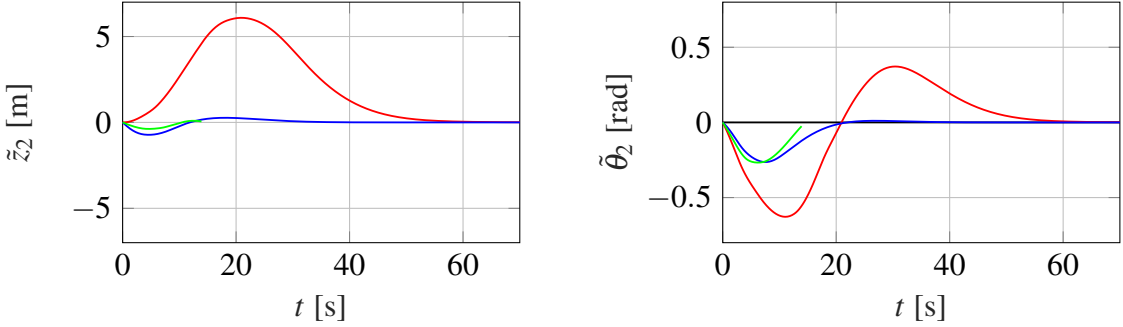}
		};
		\end{tikzpicture}
		\label{c9:fig:straight_controls}
	}
	\quad
	\setlength\figureheight{0.2\columnwidth}
	\setlength\figurewidth{0.33\columnwidth}
	\subfloat[][Lateral error of trailer 2 from $\tilde x^i_{S}$.]{
		\begin{tikzpicture}
		\node[anchor=south west] (myplot) at (0,0) {
%
%
\begin{tikzpicture}

\begin{axis}[%
width=\figurewidth,
height=\figureheight,
at={(0\figurewidth,0\figureheight)},
scale only axis,
xmin=0,
xmax=70,
xlabel={$t$ [s]},
ymin=-7,
ymax=7,
xlabel style={font=\color{white!15!black},at={(axis description cs:0.5,-0.17)},anchor=north},
ylabel style={font=\color{white!15!black},at={(axis description cs:-0.2,.5)},anchor=south},
ylabel={$\tilde z_2$ [m]},
axis background/.style={fill=white},
xmajorgrids,
ymajorgrids
]
\addplot [color=red, line width=0.7pt, forget plot]
  table[row sep=crcr]{%
0	0\\
0.5	0.00784624814450297\\
1	0.0317050999199893\\
1.5	0.0715444831367478\\
2	0.126672013043901\\
2.59999999999999	0.211080840166076\\
3.3	0.329662363512625\\
4.09999999999999	0.481951005005911\\
4.8	0.63207291071059\\
5.3	0.753969177563079\\
5.8	0.893010397823943\\
6.2	1.01952635122606\\
6.7	1.19684125977852\\
7.3	1.42780461752328\\
7.90000000000001	1.67648523302276\\
8.59999999999999	1.98686603765793\\
9.3	2.31595202431721\\
10.2	2.7598002960892\\
12.9	4.10974484886647\\
14.4	4.82369330757722\\
14.9	5.03508651011877\\
15.4	5.22714330658103\\
15.9	5.39802028754016\\
16.4	5.54673452030504\\
16.9	5.67319724193661\\
17.4	5.77811877770689\\
17.9	5.86283689934014\\
18.4	5.93039184331198\\
19	5.99400021652193\\
19.6	6.04185960102909\\
20.2	6.07440402019213\\
20.8	6.08892694016401\\
21.4	6.0836586411432\\
22	6.06010618462473\\
22.6	6.01897478469699\\
23.2	5.96073969971401\\
23.8	5.88567917327262\\
24.4	5.79371678680342\\
25	5.68476423829711\\
25.6	5.55889072833975\\
26.2	5.41645109767103\\
26.8	5.25817155262673\\
27.5	5.05504266803973\\
28.2	4.83444759455072\\
29	4.56506135247474\\
30.1	4.17442413701374\\
32.9	3.16818510934073\\
33.9	2.83184127921957\\
34.8	2.54670745276702\\
35.7	2.28003206922912\\
36.6	2.03257925931854\\
37.5	1.80447777098905\\
38.4	1.5954504203962\\
39.3	1.40496361529451\\
40.2	1.23231536770616\\
41.1	1.07668376775437\\
42	0.937156778691005\\
43	0.799818495926118\\
44	0.679726620561681\\
45.1	0.565787522474238\\
46.2	0.468952344835884\\
47.4	0.380521048927051\\
48.7	0.302229209115325\\
50.1	0.235047869757665\\
51.7	0.175963509463216\\
53.5	0.127083353339941\\
55.6	0.087387808032986\\
58.2	0.0558258936997476\\
61.7	0.0317403641285381\\
66.9	0.0148721546225232\\
70.1	0.00959329232425432\\
};
\addplot [color=blue, line width=0.7pt, forget plot]
  table[row sep=crcr]{%
0	0\\
0.700000000000003	-0.206090627869486\\
1.3	-0.363289433463095\\
1.8	-0.473136945823157\\
2.3	-0.562032219068783\\
2.8	-0.6299112752328\\
3.3	-0.677268866777069\\
3.8	-0.70590702960817\\
4.40000000000001	-0.720617027332651\\
5	-0.718104468066286\\
5.59999999999999	-0.700260737097182\\
6.2	-0.667145109650377\\
6.90000000000001	-0.610264961078627\\
7.59999999999999	-0.537335293444301\\
8.40000000000001	-0.4361680259791\\
11.2	-0.0661618793456569\\
12.1	0.0265265350419384\\
13	0.102047698028684\\
13.9	0.161431845771133\\
14.9	0.210222004936981\\
16	0.245418807309576\\
17.2	0.265055505423334\\
18.5	0.268578907003089\\
20.1	0.254459935337749\\
22.3	0.21549430844523\\
29.3	0.0788713631475986\\
32.4	0.0420127453397328\\
35.9	0.0180603803637496\\
40.6	0.00451744748227156\\
49.3	0.000266037001935615\\
70.1	5.12792753113445e-05\\
};
\addplot [color=green, line width=0.7pt, forget plot]
  table[row sep=crcr]{%
0	0\\
1.2	-0.154629052267167\\
2.1	-0.25230450281073\\
2.8	-0.312445227978227\\
3.5	-0.35471537328503\\
4.1	-0.374223609273514\\
4.8	-0.378129400757889\\
5.6	-0.365307204035206\\
6.5	-0.33446323251933\\
7.2	-0.295428983119562\\
7.8	-0.246212437266985\\
10.9	0.0399905255434181\\
11.6	0.0776293207209484\\
12.2	0.094048109444909\\
12.8	0.0938091810265647\\
13.5	0.0748714258597865\\
13.9	0.0576409385976362\\
};
\end{axis}
\end{tikzpicture}%
		};
		\end{tikzpicture}
		\label{c9:fig:z2_traj_straight}
	}
	~
	\setlength\figureheight{0.2\columnwidth}
	\setlength\figurewidth{0.33\columnwidth}
	\subfloat[][Heading error of trailer 2 from $\tilde x^i_{S}$.]{
		\begin{tikzpicture}
		\node[anchor=south west] (myplot) at (0,0) {
%
%
\begin{tikzpicture}

\begin{axis}[%
width=\figurewidth,
height=\figureheight,
at={(0\figurewidth,0\figureheight)},
scale only axis,
xmin=0,
xmax=70,
xlabel={$t$ [s]},
ymin=-0.8,
ymax=0.8,
xlabel style={font=\color{white!15!black},at={(axis description cs:0.5,-0.17)},anchor=north},
ylabel style={font=\color{white!15!black},at={(axis description cs:-0.2,.5)},anchor=south},
ylabel={$\tilde \theta_2$ [rad]},
axis background/.style={fill=white},
xmajorgrids,
ymajorgrids
]
\addplot [color=black, line width=0.7pt, forget plot]
table[row sep=crcr]{%
	0	0\\
	70	0\\
};
\addplot [color=red, line width=0.7pt, forget plot]
  table[row sep=crcr]{%
0	0\\
0.400000000000006	-0.0304346662744308\\
0.799999999999997	-0.0634467167565873\\
1.2	-0.098669021649016\\
1.7	-0.145046380529095\\
3	-0.267802000646697\\
3.3	-0.293893045529359\\
3.59999999999999	-0.318095064839184\\
3.90000000000001	-0.340484782507119\\
4.3	-0.367984284213819\\
4.7	-0.393344806164549\\
5.2	-0.422959634297413\\
5.90000000000001	-0.462240957635871\\
6.7	-0.50539179851009\\
7.2	-0.530051411812664\\
7.7	-0.552381484385251\\
8.2	-0.572207632211985\\
8.59999999999999	-0.586142912874919\\
9	-0.598260183599223\\
9.40000000000001	-0.608457301021573\\
9.8	-0.616628881665179\\
10.2	-0.62266700661209\\
10.6	-0.626462174059441\\
11	-0.62790454676437\\
11.4	-0.626885544272596\\
11.8	-0.623299830262397\\
12.2	-0.617047742000437\\
12.6	-0.60803820019899\\
12.9	-0.59936560335629\\
13.2	-0.588655110878648\\
13.5	-0.57571100339527\\
13.8	-0.560893451745926\\
14.1	-0.544318437300788\\
14.4	-0.526018207232269\\
14.7	-0.506044450917059\\
15.1	-0.476941745364869\\
15.5	-0.445249936278913\\
15.9	-0.411320564471723\\
16.4	-0.366494685051649\\
17.8	-0.238503296750338\\
18.2	-0.205227589180296\\
18.7	-0.166216878146727\\
19.4	-0.114482799098255\\
20.3	-0.0480443234970949\\
21.1	0.0123115752171827\\
21.7	0.0548833061426137\\
22.3	0.0952065811138283\\
22.9	0.13324668543116\\
23.4	0.163021738608577\\
23.9	0.191015201105529\\
24.4	0.217208909852701\\
24.9	0.241565225447104\\
25.4	0.264028223923006\\
25.9	0.284528833317566\\
26.4	0.302991942315145\\
26.9	0.319344093012248\\
27.4	0.333520642619163\\
27.9	0.345473544984898\\
28.4	0.355176891718969\\
28.9	0.362631320201572\\
29.4	0.36786694338825\\
29.9	0.370944456381793\\
30.4	0.371954324922811\\
30.9	0.371014133533535\\
31.5	0.367510939193437\\
32.1	0.361678071899476\\
32.7	0.353812184575247\\
33.4	0.34247033959771\\
34.2	0.32722720843384\\
35.1	0.307944350535379\\
36.3	0.280056532702673\\
40.2	0.187727502159859\\
41.5	0.159956332968108\\
42.7	0.136534619830456\\
43.8	0.117103550680682\\
44.9	0.0996780025857049\\
46	0.0842384456036598\\
47.2	0.0695721945870389\\
48.4	0.057031736199491\\
49.7	0.045630383606948\\
51.1	0.0356015001295447\\
52.6	0.0270762528029991\\
54.2	0.0200804264224814\\
56.1	0.0139993535455005\\
58.3	0.0092068302525945\\
61	0.00556474800619355\\
64.6	0.0029737925014075\\
70.1	0.00130375230635593\\
};
\addplot [color=blue, line width=0.7pt, forget plot]
  table[row sep=crcr]{%
0	0\\
0.400000000000006	-0.0147564689132764\\
0.799999999999997	-0.0316653171489918\\
1.2	-0.0504551098616162\\
1.59999999999999	-0.0711228464795823\\
2.2	-0.104488461841513\\
3	-0.149169420182744\\
3.40000000000001	-0.169385830089865\\
3.8	-0.187097956592936\\
4.2	-0.202490469145317\\
4.59999999999999	-0.215845575225003\\
5.09999999999999	-0.230056798235424\\
5.59999999999999	-0.241818393688519\\
6.09999999999999	-0.251278067468917\\
6.59999999999999	-0.25834576908926\\
7	-0.26192321873765\\
7.40000000000001	-0.263556652882386\\
7.8	-0.26322718099739\\
8.2	-0.260938362708373\\
8.59999999999999	-0.25680454231491\\
9.09999999999999	-0.249201144494478\\
9.60000000000001	-0.239180117505768\\
10.2	-0.224688163361819\\
11	-0.202980660606528\\
13.9	-0.122086820514554\\
14.8	-0.0997617851004122\\
15.7	-0.0795219925462618\\
16.6	-0.0615460264239545\\
17.5	-0.045902887322768\\
18.4	-0.0325703850537025\\
19.3	-0.0214511714944337\\
20.2	-0.0123907518443076\\
21.2	-0.00450205617555355\\
22.3	0.00187049673171202\\
23.5	0.00653250381053283\\
24.8	0.00945824613555146\\
26.4	0.0108549368853375\\
28.5	0.0103761627328112\\
32.3	0.0068774754874994\\
37.8	0.00232894439307074\\
42.9	0.000465789821120666\\
51.2	-7.09100521021355e-05\\
70.1	5.76658914042127e-06\\
};
\addplot [color=green, line width=0.7pt, forget plot]
  table[row sep=crcr]{%
0	0\\
0.5	-0.0291045348111361\\
1	-0.0604883558258624\\
1.7	-0.106844264462273\\
2.5	-0.159804826681869\\
2.9	-0.1845134439374\\
3.3	-0.206859513439744\\
3.6	-0.221525590284733\\
3.9	-0.234026257138401\\
4.2	-0.244123680187341\\
4.5	-0.251751583197249\\
4.8	-0.257394847215478\\
5.2	-0.262550275734551\\
5.6	-0.265613720816265\\
6.1	-0.267195068957902\\
6.6	-0.266626498204468\\
7.1	-0.263743524212886\\
7.5	-0.25931364533192\\
7.9	-0.252493667433601\\
8.4	-0.241405166675234\\
8.9	-0.227889646786027\\
9.4	-0.212201997032393\\
10	-0.190905891153632\\
10.6	-0.167410851974758\\
11.3	-0.137959354988597\\
12.5	-0.0848530080570953\\
13.4	-0.0455420734528058\\
13.9	-0.0253005850998669\\
};
\end{axis}
\end{tikzpicture}%
		};
		\end{tikzpicture}
		\label{c9:fig:th2_traj_straight}
	}	
	\caption{Results from path following of a straight nominal path ($y_{2r}=0$) in backward motion from initial state $\tilde x^i_{S}$ in Figure~\ref{c9:fig:xy_path_straight} with $(\beta_2^i,\beta_1^i)=(0.6,-0.6)$ rad using MS2T-LQ (green lines), SS2T-MPC (red lines) and  MS2T-MPC (blue lines). In Figure~\ref{c9:fig:b23_traj_straight}, gray box illustrates the joint-angle constraints.}
	\label{c9:fig:sim_straight}
\end{figure}

\begin{figure}[t!]
	\vspace{-5pt}
	\centering
	\setlength\figureheight{0.75\textwidth}
	\setlength\figurewidth{0.75\textwidth} 
	\begin{tikzpicture}
	\node[anchor=south west] (myplot) at (0,0) {
		\input{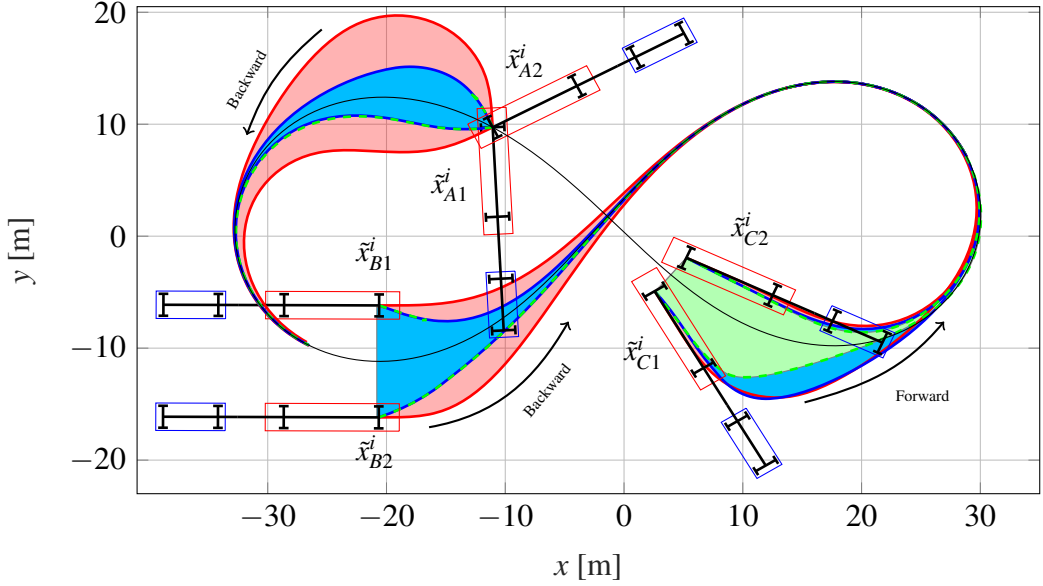}
	};
	\begin{scope}[x={(myplot.south east)}, y={(myplot.north west)}]
	\node[text=black] at (0.72,0.6) {\small $\tilde x^i_{C2}$};
	\node[text=black] at (0.62,0.4) {\small $\tilde x^i_{C1}$};
	\node[text=black] at (0.51,0.87) {\small $\tilde x^i_{A2}$};
	\node[text=black] at (0.44,0.67) {\small $\tilde x^i_{A1}$};
	\node[text=black] at (0.37,0.24) {\small $\tilde x^i_{B2}$};
	\node[text=black] at (0.37,0.56) {\small $\tilde x^i_{B1}$};
	\draw[<-, thick] (0.25,0.75) to [out=70,in=220] (0.32,0.92);  
	\node[rotate=60] at (0.25,0.84) {\tiny Backward};
	\draw[->, thick] (0.42,0.28) to [out=10,in=240] (0.55,0.45);  
	\node[rotate=50] at (0.53,0.34) {\tiny Backward};
	\draw[->, thick] (0.77,0.32) to [out=15,in=230] (0.9,0.45);  
	\node[rotate=0] at (0.88,0.33) {\tiny Forward};
	\end{scope}
	\end{tikzpicture}
	\vspace{-10pt}
	\caption{Path following of a figure-eight path in backward (Scenario A and Scenario B) and in forward motion (Scenario C) from perturbed initial states using MS2T-LQ (green lines), SS2T-MPC (red lines) and  MS2T-MPC (blue lines). In Scenario A, heading error $\tilde\theta^i_2\in[-1,\hspace{1pt}1]$ rad, in Scenario B, lateral error $\tilde z^i_2\in[-5,\hspace{1pt}5]$ m and in Scenario C, both lateral error $\tilde z^i_2\in[-2,\hspace{1pt}2]$ m and heading error $\tilde\theta^i_2\in[-0.3,\hspace{1pt}0.3]$ rad. The blue, red and green sets represent the convex envelope of the trajectories for the position of \mbox{trailer 2} using MS2T-MPC, SS2T-MPC, and MS2T-LQ, respectively. From some high-lighted initial states, MS2T-LQ leads to jack-knife (see Figure~\ref{c9:fig:b23_traj_b23_eight_rev}).}
	\label{c9:fig:xy_path_eight}
\end{figure}

The convex envelopes of the trajectories for position $(x_2(\cdotp),y_2(\cdotp))$ using MS2T-MPC and SS2T-MPC with initial joint-angle errors $\tilde\beta^i_1,\tilde\beta^i_2\in[-0.6,0.6]$ rad are illustrated in Figure~\ref{c9:fig:xy_path_straight}. The convex envelope for MS2T-LQ is not presented since the vehicle enters a jack-knife state from some initial states. From Figure~\ref{c9:fig:xy_path_straight}, it is clear that the transient response for MS2T-MPC yields a significantly smaller maximum overshoot in the lateral error of trailer 2 (\mbox{$0.26$ m}) compared to SS2T-MPC ($6.1$ m). The reason for this can be seen in Figure~\ref{c9:fig:straight_controls}, where the control-input trajectories are plotted from initial state $\tilde x^i_{S}$ with $(\tilde\beta_2^i,\tilde\beta_1^i)=(0.6,-0.6)$~rad. The results show that MS2T-MPC uses a positive trailer-steering angle $\gamma_2$ (blue dashed line) to compensate for the initially positive value of $\beta_2$ (see Figure~\ref{c9:fig:b23_traj_straight}). Moreover, as can be seen Figure~\ref{c9:fig:th2_traj_straight}, the trailer-steering angle is also used by MS2T-MPC to reduce the maximum overshoot in the heading error of trailer 2 \mbox{($0.26$ rad)} compared to SS2T-MPC \mbox{($0.62$ rad)}. 
Additionally, the joint-angle trajectories from initial state $\tilde x^i_{S}$ are plotted in Figure~\ref{c9:fig:b23_traj_straight}. As can be seen, MS2T-LQ is not able to stabilize the vehicle due to the input constraints and jackknifing occurs almost instantly. As a comparison, both MPC controllers are able to make the system converge to the straight nominal path, but SS2T-MPC has to initially violate the soft joint-angle constraints, which is not the case for MS2T-MPC.

\begin{figure}[t!]
	\centering
	\captionsetup[subfloat]{captionskip=-2pt} 
	\setlength\figureheight{0.22\textwidth}
	\setlength\figurewidth{0.33\textwidth}
	\subfloat[][Joint-angle trajectories from initial state $\tilde x^i_{A1}$ in Figure~\ref{c9:fig:xy_path_eight}. Initial state denoted by black star and nominal path  $(\beta_{1r}(\cdotp),\beta_{2r}(\cdotp))$ by black line.]{
		\begin{tikzpicture}
		\node[anchor=south west] (myplot) at (0,0) {
			\input{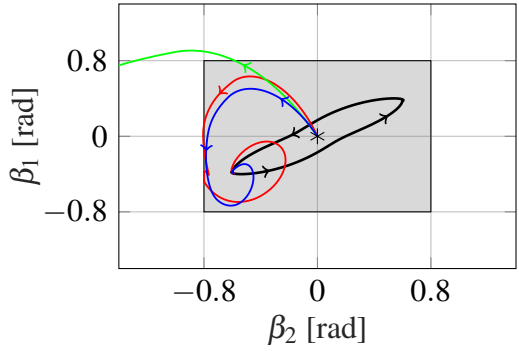}
		};
		\end{tikzpicture}
		\label{c9:fig:b23_traj_b23_eight_rev}
	}
	~
	\setlength\figureheight{0.22\textwidth}
	\setlength\figurewidth{0.33\textwidth}
	\subfloat[][Joint-angle trajectories from initial state $\tilde x^i_{C1}$ in Figure~\ref{c9:fig:xy_path_eight}. Initial state denoted by black star and nominal path in $(\beta_{1r}(\cdotp),\beta_{2r}(\cdotp))$ by black line.]{
		\begin{tikzpicture}
		\node[anchor=south west] (myplot) at (0,0) {
			\input{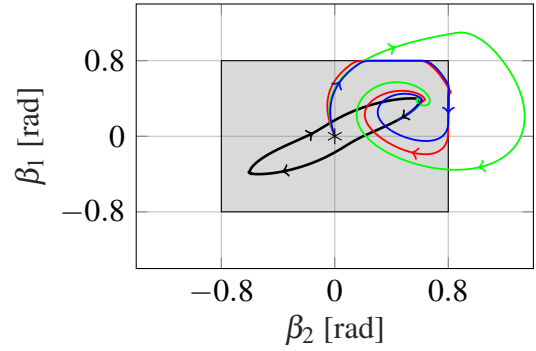}
		};
		\end{tikzpicture}
		\label{c9:fig:b23_traj_b23_eight_fwd}
	}
	\quad
	\setlength\figureheight{0.16\textwidth}
	\setlength\figurewidth{0.33\textwidth}
	\subfloat[][Heading error of trailer 2 from $\tilde x^i_{A1}$.]{
		\begin{tikzpicture}
		\node[anchor=south west] (myplot) at (0,0) {
			\input{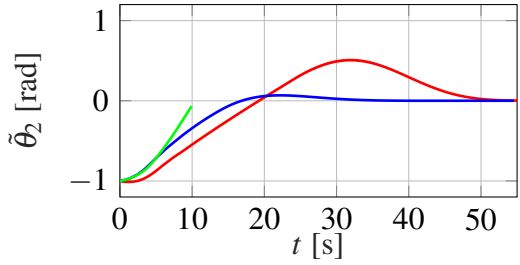}
		};
		\end{tikzpicture}
		\label{c9:fig:z2_traj_eight_rev}
	}
	~
	\subfloat[][Lateral error of trailer 2 from $\tilde x^i_{A1}$.]{
		\begin{tikzpicture}
		\node[anchor=south west] (myplot) at (0,0) {
%
%
%
\begin{tikzpicture}

\begin{axis}[%
width=\figurewidth,
height=\figureheight,
at={(0\figurewidth,0\figureheight)},
scale only axis,
xmin=0,
xmax=55,
xlabel={$t$ [s]},
ymin=-9,
ymax=9,
xlabel style={font=\color{white!15!black},at={(axis description cs:0.5,-0.13)},anchor=north},
ylabel style={font=\color{white!15!black},at={(axis description cs:-0.17,.5)},anchor=south},
ylabel={$\tilde z_2$ [m]},
axis background/.style={fill=white},
xmajorgrids,
ymajorgrids
]
\addplot [color=red, line width=1.0pt, forget plot]
  table[row sep=crcr]{%
0	6.86717796867242e-05\\
0.619999999999997	0.501024110568252\\
1.26	0.973358962855599\\
2.14	1.57866583448563\\
5.04	3.51439042815146\\
5.71	3.89818723043592\\
6.36	4.2334375943737\\
7.03	4.5423800565312\\
8	4.95133907061347\\
9.69	5.62911185012466\\
10.46	5.89321792903973\\
11.47	6.20536169902648\\
12.45	6.47344076599949\\
13.39	6.69755381884436\\
14.33	6.88887250552558\\
15.26	7.04540243399293\\
16.15	7.1646472496118\\
17.06	7.25547779271609\\
18.02	7.31732294896916\\
18.97	7.3442457415352\\
19.9	7.33770343843184\\
20.8	7.30069144497993\\
21.74	7.22998703607804\\
22.65	7.13076389320665\\
23.62	6.99202066856914\\
24.64	6.81216223700209\\
25.7	6.59332195949536\\
26.73	6.34586217464616\\
27.72	6.07208681708395\\
28.7	5.7643039457668\\
29.63	5.43866544481513\\
30.69	5.03130105580753\\
31.79	4.5748960029866\\
33.53	3.81143336903546\\
35.71	2.85702236273151\\
36.93	2.35866258122496\\
37.98	1.96403702958204\\
38.93	1.63972136924848\\
39.86	1.35517813733943\\
40.81	1.09957349147987\\
41.76	0.879514763809098\\
42.68	0.699258792622551\\
43.62	0.546554545418694\\
44.66	0.411488168996939\\
45.71	0.306781080496371\\
46.93	0.21855807556372\\
48.23	0.155607468846448\\
49.93	0.107563723368791\\
52.34	0.0783289658543396\\
55.01	0.068309561486501\\
};
\addplot [color=blue, line width=1.0pt, forget plot]
  table[row sep=crcr]{%
0	6.86717796867242e-05\\
0.630000000000003	0.387460107989142\\
1.31	0.767449687552251\\
2.16	1.20385150468805\\
3.11	1.65664387751264\\
3.86	1.98152693104272\\
4.42	2.18900130244529\\
4.99	2.36601132627323\\
5.58	2.51481077728975\\
6.23	2.6451952896933\\
6.98	2.76185775191897\\
7.79	2.85529497040505\\
8.63	2.91998259943979\\
9.44	2.95184631709806\\
10.28	2.9528547825195\\
11.12	2.92123159087924\\
11.95	2.85785097336307\\
12.77	2.76349506346137\\
13.57	2.64098551798852\\
14.66	2.43611844590354\\
17.23	1.91635621717133\\
19.62	1.44246393749923\\
21.35	1.13437407547816\\
22.74	0.91757804146102\\
24.12	0.732499507463267\\
25.5	0.577703961953887\\
27.01	0.440984561351847\\
28.64	0.327282762480266\\
30.43	0.23583835706593\\
32.59	0.160876465098994\\
35.04	0.107784968397205\\
38.35	0.0674177820402448\\
44.23	0.0337156690156419\\
54.54	0.00966959881937868\\
};
\addplot [color=green, line width=1.0pt, forget plot]
  table[row sep=crcr]{%
0	6.86717796885006e-05\\
0.6	0.368067124928112\\
1.07	0.619718003803079\\
1.58	0.855974326812129\\
2.16	1.08794056915085\\
2.81	1.31297226480056\\
3.54	1.53158656097711\\
4.22	1.70261901359976\\
4.73	1.8008136766213\\
5.26	1.86936149233115\\
5.83	1.90986022711144\\
6.53	1.92642842514258\\
7.26	1.91172701029358\\
8.04	1.8637411249622\\
9	1.77092199264093\\
9.93	1.6683510776507\\
};
\end{axis}
\end{tikzpicture}%
		};
		\end{tikzpicture}
		\label{c9:fig:th2_traj_eight_rev}
	}
	\quad
	\setlength\figureheight{0.14\textwidth}
	\setlength\figurewidth{0.7\textwidth} 
	\subfloat[][Control inputs $\kappa_0$ (solid) and $\gamma_2$ (dashed) from $\tilde x^i_{A1}$, and their limits $\bar\kappa_0$ (dashed-dotted black) and $\bar\gamma_2$ (dotted black).]{
		\begin{tikzpicture}
		\node[anchor=south west] (myplot) at (0,0) {
			\input{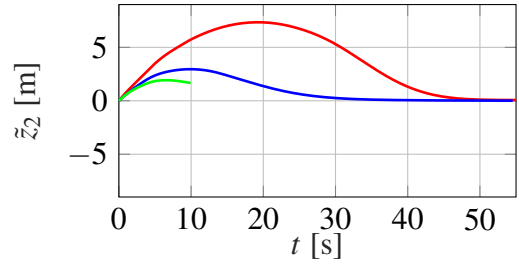}
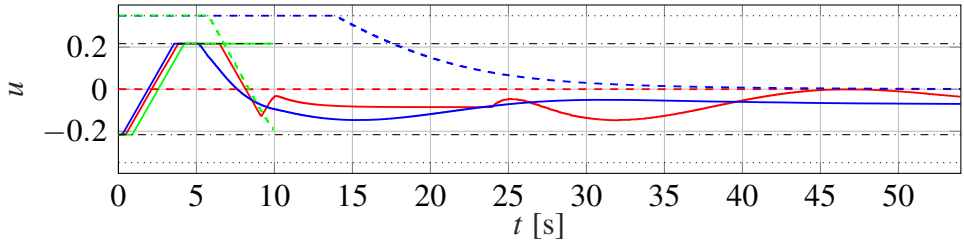
		};
		\end{tikzpicture}
		\label{c9:fig:control_eight_rev}
	}	
	\caption{Results from path following of a figure-eight path in backward (Scenario A and Scenario B) and in forward motion (Scenario C) from perturbed initial states using MS2T-LQ (green lines), SS2T-MPC (red lines) and  MS2T-MPC (blue lines). In Scenario A (see Figure~\ref{c9:fig:xy_path_eight}), heading error $\tilde\theta^i_2\in[-1,\hspace{1pt}1]$ rad, in Scenario B, lateral error $\tilde z^i_2\in[-5,\hspace{1pt}5]$ m and in Scenario C, both lateral error $\tilde z^i_2\in[-2,\hspace{1pt}2]$ m and heading error $\tilde\theta^i_2\in[-0.3,\hspace{1pt}0.3]$ rad. In Figure~\ref{c9:fig:b23_traj_b23_eight_rev}-\ref{c9:fig:b23_traj_b23_eight_fwd}, the gray box illustrates the used joint-angle constraints.}
	\label{c9:fig:path_eight}
\end{figure}

The second set of simulations involve backward tracking ($v_0=-1$ m/s) and forward tracking ($v_0=1$ m/s) of a figure-eight nominal path in ($x_{2r}(\cdotp),y_{2r}(\cdotp)$), which has been computed as described in~\cite{LjungqvistJFR2019}. Also in this set of simulations, the initial state $\tilde x(0)$ is perturbed to compare the performance of the controllers. The simulation results are presented in Figure~\ref{c9:fig:xy_path_eight}-\ref{c9:fig:path_eight}. Scenario~A involves a heading error $\tilde \theta^i_2\in[-1,1]$ rad and Scenario~B a lateral error $\tilde z^i_2\in[-5,5]$ m, both in backward motion. The same conclusions as for the first set of simulations can be drawn, the MS2T-MPC uses the trailer-steering angle to reduce the overshoot and convergence times for $\tilde z_2$ and $\tilde \theta_2$ compared to SS2T-MPC (see, e.g., Figure~\ref{c9:fig:z2_traj_eight_rev}--\ref{c9:fig:th2_traj_eight_rev}). Moreover, MS2T-LQ is not able to stabilize the system and jack-knife occurs almost instantly from some initial states, e.g., $\tilde x^i_{A1}$ with $\tilde \theta^i_2=-1$ rad (see Figure~\ref{c9:fig:b23_traj_b23_eight_rev}) and $\tilde x^i_{B1}$ with $\tilde z^i_2=-5$ m. Finally, Scenario C involves both a lateral $\tilde z^i_2\in[-2,2]$ m and an initial heading error $\tilde \theta^i_2\in[-0.3,0.3]$ rad in forward motion. As can be seen in Figure~\ref{c9:fig:xy_path_eight}, the convex envelope of the trajectories $(x_2(\cdotp),y_2(\cdotp))$ is smallest for MS2T-LQ (green set). However, the joint-angle trajectories (see Figure~\ref{c9:fig:b23_traj_b23_eight_fwd}) in the MS2T-LQ case are drastically violating their constraints at some parts of the maneuvers, which is neither the case for MS2T-MPC nor SS2T-MPC. As a final note, the average computation time in Gurobi for the proposed MS2T-MPC is 35 ms compared to 25 ms for SS2T-MPC which is less than the sampling time $T_s=100$ ms of the controllers.

\section{Conclusions}
\label{c9:sec:conclusions}
A model predictive path-following controller is proposed for multi-steered articulated vehicles composed of a car-like tractor and an arbitrary number of off/on-axle hitched trailers with steerable/non-steerable wheels. The proposed MPC controller uses a path-following error model of the vehicle for predictions, is designed to satisfy physically constraints on states and control inputs, and is tailored to follow nominal paths that contain full state and control-input information. The performance of the proposed path-following controller is evaluated in a set of practically relevant scenarios for a multi-steered 2-trailer with a car-like tractor where the last trailer is steerable. In simulations, it is shown that the proposed controller outperforms a linear quadratic controller and efficiently exploits the additional trailer-steering capability, while recovering from non-trivial initial states in backward motion.   

As future work, we would like develop a motion planner and evaluate the framework in real-work experiments on a full-scale test vehicle.  

\section*{Appendix A}
We start by deriving the functions $f_{\beta_1}$, $f_{\beta_2}$ and $\kappa_2$ describing the path-following error model for the specific MS2T with a car-like tractor~\eqref{c9:eq:model_frenet_frame_2Ts}. The matrices $J_1$ and $J_2$ describing the longitudinal and angular velocity transformations between neighboring vehicle segments~\eqref{c9:eq:velocity_transformation} are 
\begin{align}
\label{c9:MS2T_vel_matrices}
\begin{split}
J_1(\beta_1,0,0)=\begin{bmatrix}
-\frac{M_{1}}{L_{1}}\cos\beta_{1} & \frac{\sin\beta_1}{L_{1}} \\[10pt]
M_{1}\sin\beta_1 & \cos\beta_{1}
\end{bmatrix}, \quad J_2(\beta_2,\gamma_2,0) = \begin{bmatrix}
0 & \frac{\sin{(\beta_{2}-\gamma_{2})}}{L_{2}\cos\gamma_{2}} \\[10pt]
0 & \frac{\cos\beta_{2}}{\cos\gamma_{2}}
\end{bmatrix},
\end{split}
\end{align}
since $M_2=\gamma_0=\gamma_1=0$. Thus, the velocity transformation from trailer 2 to the car-like tractor~\eqref{c9:eq:velocity_0N_trans} is  
\begin{align}
\label{c9:MS2T_vel_trans}
f_{v_2}(\beta_1,\beta_2,u) &= d^T J_2 J_1\begin{bmatrix}
\kappa_0 \\ 1 
\end{bmatrix} = \frac{\cos\beta_2}{\cos\gamma_2}\left(M_1\sin\beta_1\kappa_0 + \cos\beta_1\right), 
\end{align}
and the curvature of trailer 2~\eqref{c9:eq:curvature_Ntrailer} is
\begin{align}
\kappa_2(\beta_2,\gamma_2) &= \frac{c^T  J_2 J_1\begin{bmatrix}
	\kappa_0 \\ 1 
	\end{bmatrix}}{d^T  J_2  J_1\begin{bmatrix}
	\kappa_0 \\ 1 
	\end{bmatrix}} = \frac{\sin(\beta_2-\gamma_2)}{L_2\cos\beta_2}.
\label{c9:curvature-2trailer}
\end{align} 
Finally, using~\eqref{c9:MS2T_vel_matrices}--\eqref{c9:curvature-2trailer} the functions describing the joint-angle kinematics~\eqref{c9:eq:joint_angle_kinematics} are 
\begin{align}
f_{\beta_2}(\beta_1,\beta_2,u) =&\frac{\cos\gamma_2\left(\frac{\sin\beta_1}{L_1} - \frac{M_1}{L_1}\cos\beta_1\kappa_0 \right)}{\cos\beta_2(M_1\sin\beta_1\kappa_0+\cos\beta_1)} -  \frac{\sin(\beta_2-\gamma_2)}{L_2\cos\beta_2}, \nonumber \\
f_{\beta_1}(\beta_1,\beta_2,u) =&\frac{\cos\gamma_2\left(\kappa_0 -\frac{\sin\beta_1}{L_1} + \frac{M_1}{L_1}\cos\beta_1\kappa_0 \right)}{\cos\beta_2(M_1\sin\beta_1\kappa_0+\cos\beta_1)}.
\label{c9:joint-angle-kinematics-2trailer}
\end{align}
The Jacobian linearization of the nonlinear path-following error model~\eqref{c9:eq:model_frenet_frame_2Ts} around the origin $(\tilde x,\tilde u) = (0,0)$ can be represented as in~\eqref{c9:linear_cont_model}, where the matrices $A(s)$ and $B(s)$ have the following structure
\begin{align}\label{c9:eq:A}
A(s) = \frac{\partial \tilde f(s,0,0)}{\partial \tilde x} = \begin{bmatrix}
0 & 1 & 0 & 0 \\
a_{21}(s) & 0 & a_{23}(s) & 0 \\
a_{31}(s) & 0 & a_{33}(s) & a_{34}(s) \\
a_{41}(s)& 0 & a_{43}(s) & a_{44}(s) 
\end{bmatrix}, 
\end{align}
and 
\begin{align}\label{c9:eq:B}
B(s) = \frac{\partial \tilde f(s,0,0)}{\partial \tilde u} = \begin{bmatrix}
0 & 1 \\ 0 & b_{22}(s) \\ b_{31}(s) & b_{32}(s)\\ b_{41}(s) & b_{42}(s)
\end{bmatrix}, 
\end{align}
where
\begingroup\makeatletter\def\f@size{9.2}\check@mathfonts
\begin{subequations}
	\begin{align*}
	a_{21}(s) =& -\frac{\sin^2{(\beta_{2r}-\gamma_{2r})}}{L_2^2\cos^2{\beta_{2r}}}, \\
	a_{23}(s) =& \frac{\cos(\beta_{2r}-\gamma_{2r})}{L_2\cos\beta_{2r}} + \frac{\sin(\beta_{2r}-\gamma_{2r})\tan\beta_{2r}}{L_2\cos\beta_{2r}}, \\
	a_{31}(s) =& -\frac{\sin(\beta_{2r}-\gamma_{2r})}{L_2\cos^2\beta_{2r}}\left(
	\frac{\cos\gamma_{2r}(\sin\beta_{1r}-\kappa_{0r}M_1\cos\beta_{1r})}{L_1(\cos\beta_{1r} + \kappa_{0r}M_1\sin{\beta_{1r}})} - \frac{\sin(\beta_{2r}-\gamma_{2r})}{L_2} \right),  \\
	a_{33}(s) =& \cos\gamma_{2r}\left(\frac{\sin\beta_{2r}(\sin\beta_{1r}-\kappa_{0r} M_1 \cos\beta_{1r})}{L_1 \cos\beta^2_{2r} (\cos\beta_{1r} + \kappa_{0r} M_1 \sin\beta_{1r})} -\frac{1}{\cos^2\beta_{2r}L_2}\right), \\
	a_{34}(s) =& \frac{\cos\gamma_{2r}(1+\kappa_{0r}^2M_1^2)}{L_1\cos\beta_{2r}(\cos\beta_{1r}+\kappa_{0r}M_1\sin{\beta_{1r}})^2}, \\
	a_{41}(s) =& 
	-\frac{\cos\gamma_{2r}(\kappa_{0r}L_1 - \sin\beta_{1r} + M_1\cos\beta_{1r}\kappa_{0r})}{L_1 L_2 \cos^2\beta_{2r}(\cos\beta_{1r} + M_1\kappa_{0r}\sin\beta_{1r})}, \\
	a_{43}(s) =& \frac{\cos\gamma_{2r}\tan\beta_{2r}}{L_1}\left(\frac{\kappa_{0r}L_1 + \kappa_{0r} M_1 \cos\beta_{1r} - \sin\beta_{1r}}{\cos\beta_{2r}(\cos\beta_{1r} +\kappa_{0r}M_1\sin\beta_{1r})}\right), \\
	a_{44}(s) =& \frac{1 + \kappa_{0r}^2 M_1^2 + \kappa_{0r}^2L_1M_1\cos\beta_{1r} -\kappa_{0r} L_1 \sin\beta_{1r}}{\sec\gamma_{2r}L_1\cos\beta_{2r}(\cos\beta_{1r}+\kappa_{0r} M_1\sin\beta_{1r})^2}, 
	\end{align*}
\end{subequations}
\endgroup
and
\begingroup\makeatletter\def\f@size{9.2}\check@mathfonts
\begin{subequations}
	\begin{align*}
	b_{22}(s) &=-\frac{\cos(\beta_{2r}-\gamma_{2r})}{L_2\cos\beta_{2r}}, \\
	b_{31}(s) &= - \frac{M_1\cos\gamma_{2r}}{L_1\cos\beta_{2r}(\cos\beta_{1r}+\kappa_{0r} M_1 \sin\beta_{1r})^2},\\
	b_{32}(s) &= \frac{\cos(\beta_{2r}-\gamma_{2r})}{L_2\cos\beta_{2r}} + \frac{(\kappa_{0r}M_1\cos\beta_{1r}-\sin\beta_{1r})\sin\gamma_{2r}}{\cos\beta_{2r}L_1(\cos\beta_{1r}+\kappa_{0r}M_1\sin\beta_{1r})},\\
	b_{41}(s) &= \frac{\cos\gamma_{2r}(M_1 + L_1\cos\beta_{1r})}{\cos\beta_{2r}L_1(\cos\beta_{1r}+\kappa_{0r} M_1 \sin\beta_{1r})^2},\\
	b_{42}(s) &= -\frac{(\kappa_{0r}L_1+\kappa_{0r}M_1\cos\beta_{1r}-\sin\beta_{1r})\sin\gamma_{2r}}{\cos\beta_{2r}L_1(\cos\beta_{1r}+\kappa_{0r}M_1\sin\beta_{1r})}.
	\end{align*}
\end{subequations}
\endgroup
Around a straight nominal path, the system matrices in~\eqref{c9:eq:A} and \eqref{c9:eq:B} simplify to
\setcounter{equation}{40}
\begin{align}\label{c9:eq:AB_straight}
 A = \begin{bmatrix}
0 & 1 & 0 & 0 \\
0 & 0 & \frac{1}{L_2} & 0 \\[6pt]
0 & 0 & -\frac{1}{L_2} & \frac{1}{L_1} \\[6pt]
0 & 0 & 0 & \frac{1}{L_1} 
\end{bmatrix}, \quad  B = \begin{bmatrix}
0 & 1 \\ 0 & -\frac{1}{L_2} \\[4pt] -\frac{M_1}{L_1} & \frac{1}{L_2}\\[6pt] \frac{M_1+L_1}{L_1} & 0
\end{bmatrix}.
\end{align}

\bibliography{root}
\bibliographystyle{abbrv}

\end{document}